\documentclass{article}
\usepackage[utf8]{inputenc}
\usepackage{mathrsfs}
\usepackage{amssymb}
\usepackage{amsmath}
\usepackage{bm}
\usepackage{cite}
\usepackage[algo2e]{algorithm2e}
\usepackage{tikz,pgfplots}
\pgfplotsset{compat=1.18}
\usepackage{amsfonts}   
\usepackage{subfigure}
\usepackage{graphicx}       % graphics
\usepackage{braket}
\usepackage{booktabs}
\usepackage{multirow}
\usepackage[margin=1in]{geometry}

\usepackage[normalem]{ulem}
\usepackage[colorlinks,linkcolor=blue,anchorcolor=blue,citecolor=red]{hyperref}
\usepackage{enumerate}
\usepackage{cases}
\usepackage{mathrsfs}
\usepackage{calligra}
\newtheorem{remark}{Remark}

\newcommand\bx{\boldsymbol{x}}

\newcommand\bu{\boldsymbol{u}}
\newcommand\bv{{\boldsymbol{v}}}
\newcommand\tbx{\tilde{\boldsymbol{x}}}
\newcommand\tbv{\tilde{\boldsymbol{v}}}

\newcommand{\Leq}{L_{\rm Eq}}
\newcommand{\Lbc}{L_{\rm BC}}
\newcommand{\Lic}{L_{\rm IC}}

\newcommand\bg{\boldsymbol{g}}
\newcommand\bs{\boldsymbol{s}}
\newcommand\bh{\boldsymbol{h}}
\newcommand\bbR{\mathbb{R}}

\newcommand\Kn{{\rm Kn}}

\newcommand\mQ{\mathcal{Q}}

\newcommand\pd[2]{\frac{\partial {#1}}{\partial {#2}}}

\newcommand{\bz}{{\bm 0}}

\newcommand{\dd}{\; \mathrm{d}}

\graphicspath{{media/}}     % organize your images and other figures under media/ folder

\numberwithin{equation}{section}

\title{Neural network representation of microflows with BGK model}

\author{
  ~~Pei Zhang\footnote{Beijing Computational Science Research Center,
  Beijing, China, email:
  \texttt{zhangpei@csrc.ac.cn}.},
  ~~Yanli Wang\footnote{Beijing Computational Science Research Center,
  Beijing, China, email:
  \texttt{ylwang@csrc.ac.cn}.}
  }

\begin{document}
\maketitle
%\tableofcontents
%\clearpage
\begin{abstract}
We consider the neural representation to solve the Boltzmann-BGK equation, especially focusing on the application in microscopic flow problems. A new dimension reduction model of the BGK equation with the flexible auxiliary distribution functions is first deduced to reduce the problem dimension. Then, a network-based ansatz that can approximate the dimension-reduced distribution with extremely high efficiency is proposed. Precisely, fully connected neural networks are utilized to avoid discretization in space and time. A specially designed loss function is employed to deal with the complex Maxwell boundary in microscopic flow problems. Moreover, strategies such as multi-scale input and Maxwellian splitting are applied to enhance the approximation efficiency further. Several classical numerical experiments, including 1D  Couette flow and Fourier flow problems and 2D duct flow and in-out flow problems are studied to demonstrate the effectiveness of this neural representation method. 

    
\end{abstract}

\section{Introduction}
\label{sec:intro}
In the fields of aerospace and micro-electro-mechanical systems (MEMS), the simulation of kinetic theory has garnered significant attention. A key area of focus within this domain is the rarefied gas dynamics \cite{shen2006rarefied}, which focuses on the gas flows where the mean free path of the gas molecules is comparable to the characteristic length scale of the system. In this case, traditional continuum fluid models including Euler and Navier–Stokes equations are no longer applicable.  Instead, the Boltzmann equation, which utilizes a distribution function to describe the behavior of dilute gases at the mesoscopic scale, provides an essential framework for studying rarefied gas dynamics. However, the high dimension and the complex quadratic collision term of the Boltzmann equation pose substantial challenges to achieving efficient and accurate numerical simulations.

For the traditional numerical methods to solve the Boltzmann equation, it can be divided into deterministic methods and stochastic methods. The direct simulation Monte Carlo (DSMC) \cite{bird1994molecular, oran1998direct, ganjaei2009new, liu2021unified} is a stochastic method, which is effective for simulating high-speed rarefied gas flows. However, they are limited by inefficiency and statistical noise when dealing with low-speed flows. On the other hand, deterministic methods offer a range of solutions that tend to be more suitable for continuous, low-speed flows, where statistical variations are less dominant. The discrete velocity method (DVM) \cite{broadwell1964study, buet1996discrete, liu2020unified} is a classical deterministic method, which discretizes the distribution function at a series of microscopic velocity points. The Fourier spectral method \cite{mouhot2006fast, wu2013deterministic, gamba2017fast} and the Hermite spectral method \cite{wang2019approximation, li2022hermite, li2023hermite} utilize global basis functions to approximate the distribution function, which leads to a higher order of convergence. Additionally, the moment method \cite{grad1949kinetic}, which was proposed by Grad, simplifies the Boltzmann equation by approximating the distribution function as a series of moments.

In recent years, advances in computing power have spurred the development of numerous neural network-based methods to solve the Boltzmann equation. These methods can generally be categorized into three types. The first involves using neural networks to create surrogate models that approximate the collision operator \cite{holloway2021acceleration, xiao2021using, miller2022neural}, significantly accelerating the computation of the collision term. The second type utilizes neural networks to learn closed-form reduced models of the Boltzmann equation \cite{han2019uniformly, huang2023machine, schotthofer2022neural, li2023learning}. Finally, within the framework of the physics-informed neural network (PINNs) \cite{raissi2019physics}, the Boltzmann equation is solved by incorporating the partial differential equation along with its initial and boundary conditions into the neural network loss function, transforming the problem into an optimization problem. Then, it is extended to solve the phonon Boltzmann transport equation integrated with the Monte Carlo method in \cite{lin2024monte}. In \cite{lou2021physicsinformed}, the PINN method was initially applied to solve the BGK model by decomposing the distribution function into equilibrium and non-equilibrium parts, each approximated by separate networks, effectively addressing both forward and inverse problems in continuum and rarefied flows. Further developments include the introduction of asymptotic-preserving neural networks (APNNs) \cite{jin2023asymptotic, jin2024asymptotic} to handle multiscale time-dependent kinetic problems, and the method of a neural sparse representation method \cite{li2024solving} to tackle both the BGK and full nonlinear collision models. Additionally, the development of the PINN-DVM model has enabled the simulation of rarefied gas flows by solving the linearized steady-state BGK equation \cite{zhang2023simulation}.

In this work, a novel training-based method termed the dimension-reduced neural representation method (DRNR) is proposed to solve the Boltzmann-BGK equation using the neural network. First, building on the foundational work in \cite{chu1965kinetictheoretic, yang1995rarefied}, we propose an enhanced new dimension reduction model by incorporating the flexible auxiliary distribution functions related to the microscopic velocity. This enhancement allows for the alignment of the microscopic velocity space dimensions of the BGK equation with the spatial space dimension when the distribution function is plane-symmetric in the related direction in the microscopic velocity space. In this framework, the reduced Maxwell boundary condition \cite{maxwell1878iii, struchtrup2005macroscopic} applied in the microscopic flow problem is also deduced. With this dimension reduction model, the dimension of the independent variables decreases and the total computational cost will be reduced greatly, although the number of the unknown distribution functions is increased. 

Then the semi-discrete dimension-reduced model of the BGK equation is deduced in the framework of the discrete velocity method (DVM) with the network-based representation proposed for the dimension-reduced distribution function, which is an efficient approximation. A fully connected neural network is utilized to parameterize the dimension-reduced distribution functions with the strategies as the multi-scale inputs and Maxwellian splitting \cite{li2024solving, jin2024asymptotic} adopted. This ensures that the network's output conforms to the multi-scale nature of the Boltzmann equation and the structural characteristics of the distribution function, thereby improving the efficiency of the approximation. A specially designed loss function with a particular focus on the reduced Maxwell boundary condition is proposed for the microscopic flow problem with an adaptive weighting strategy incorporated into the training process. A series of numerical experiments are conducted to validate the effectiveness of the proposed DRNR method, focusing primarily on the microscopic flow problem with the Maxwell boundary condition. The numerical tests include the one-dimension Couette flow and Fourier heat transfer flow as well as the two-dimension rectangular duct flows and in-out flows.

The rest of this paper is organized as follows: Sec. \ref{sec:Boltz} introduces several basic properties of the Boltzmann-BGK equation as well as the Maxwell boundary condition utilized in the microscopic flow problem and the dimension reduction method with the flexible auxiliary distribution functions. 
The detailed exposition of the DRNR method including the semi-discrete DVM model, the architecture of the network, and the special design of the loss function is provided in Sec. \ref{sec:net}, with the numerical experiments presented in Sec. \ref{sec:num}. Finally, some concluding remarks are given in Sec. \ref{sec:conclusion}.

\section{Boltzmann-BGK equation}
\label{sec:Boltz}
In this section,  several basic properties of the Boltzmann equation, the Maxwell solid wall boundary, and the dimension reduction method for the BGK-type collision model will be introduced.

\subsection{Preliminary}
\label{sec:pre}
The Boltzmann equation describes the movement of the particles from the statistical view, where $f(t, \bx, \bv)$ is the density function of the particles, with the time $t$, the spatial variable $\bx$, and the microscopic velocity $\bv$. The specific form of the Boltzmann is as below 
\begin{equation}
    \label{eq:Boltzmann}
    \pd{f}{t} + \bv \cdot \nabla_{\bx} f  = \mQ(f, f),  \qquad t \geqslant 0, \quad \bx \in \bbR^{3}, \quad \bv \in \bbR^{3},
\end{equation}
where $\mQ(f, f)$ is the collision operator with a quadratic form \cite{bird1994molecular}, defined as follows:
\begin{equation}
    \label{eq:quic_coll}
    \mathcal{Q}(f, f)=\int_{\bbR^{3}} \int_{\mathbb{S}^{2}} B\left(\bv-\bv_{{\ast}}, \sigma\right)\left[f\left(\bv_{{\ast}}^{\prime}\right) f\left(\bv^{\prime}\right)-f\left(\bv_{{\ast}}\right) f(\bv)\right] \dd \sigma \dd \boldsymbol{v}_{{\ast}},
\end{equation}
where $\bv$ and $\bv_{\ast}$ denote the velocities of two particles before the collision, while $\bv^{\prime}$ and $\bv_{\ast}^{\prime}$ represent their velocities post-collision. The post-collision velocities are determined by the conservation laws of momentum and energy during the collision process. $\sigma$ is the unit normal vector to the 2-dimensional unit sphere $\mathbb{S}^{2}$. The collision kernel $B(\cdot,\cdot )$ is a non-negative function that describes the interactions between molecules, and it is determined by the mutual potential between the particles.

Due to the complexity of the quadratic collision model, several simplified collision models have been developed, one of which is the BGK collision model \cite{bhatnagar1954model} defined as
\begin{equation}
    \label{eq:bgk_coll}
    \mathcal{Q}^{\text{BGK}}[f] = \frac{1}{\Kn}(\mathcal{M}[f]-f),
\end{equation}
where $\Kn$ is the Knudsen number \cite{struchtrup2005macroscopic,hu2017asymptotic}, describing the rarefaction of a gas and $\mathcal{M}[f]$ is a local Maxwellian distribution function, dependent on the local gas density $\rho$, macroscopic velocity $\boldsymbol{u}$, and temperature $T$, which has the following form
\begin{equation}
\label{eq:Maxwell}
    \mathcal{M}[f] = \frac{\rho}{\sqrt{2\pi T}^{3}}\exp\left(-\frac{|\bv-\boldsymbol{u}|^2}{2T}\right).
\end{equation}
The macroscopic variables such as the density $\rho$, macroscopic velocity $\boldsymbol{u}=(u_1,u_2,u_{3})^T$, and temperature $T$ could be derived from the distribution function as
\begin{equation}
    \label{eq:mac_var}
    \begin{aligned}
        &\text{Density:} &\quad&\rho(t,\bx) = \int_{\bbR^{3}} f(t,\bx,\bv) \dd\bv,\\
        &\text{Momentum:} &\quad&\boldsymbol{m}(t,\bx) \triangleq\rho(t,\bx)\boldsymbol{u}(t,\bx) = \int_{\bbR^{3}} \bv f(t,\bx,\bv) \dd\bv,\\
        &\text{Energy:} &\quad&E(t,\bx) \triangleq \frac{3}{2}\rho T + \frac{1}{2}\rho |\bu|^2= \frac{1}{2}\int_{\bbR^{3}} |\bv|^2 f(t,\bx,\bv) \dd\bv.\\
    \end{aligned}
\end{equation}
The stress tensor $\sigma_{ij}$ and heat flux $q_i$ can also be derived from the distribution function as
\begin{equation}
\label{eq:sigma_q}
    \begin{aligned}
        \sigma_{ij}(t,\bx) &= \int_{\bbR^{3}} \left((v_{i}-u_i)(v_{j}-u_j)-\frac{1}{3}\delta_{ij}|\bv-\boldsymbol{u}|^2\right)f(t,\bx,\bv) \dd\bv,\qquad &i,j = 1,2, 3,\\
        q_{i}(t,\bx) &= \frac{1}{2}\int_{\bbR^{3}} |\bv-\boldsymbol{u}|^2(v_{i}-u_i)f(t,\bx,\bv) \dd\bv, \qquad &i = 1, 2, 3.\\
    \end{aligned}
\end{equation}
So far, we have introduced the Boltzmann equation and discussed its fundamental properties. Next, we will discuss particle-wall interactions and Maxwell boundary conditions \cite{maxwell1878iii} for the Boltzmann equation.

\subsection{Maxwell boundary condition}
\label{sec:boun}
For the Boltzmann equation, several applications involve the interaction between gases and solid walls, especially in microscopic flow. For such gas-solid wall interactions, the Maxwell boundary condition \cite{maxwell1878iii, struchtrup2005macroscopic}, which combines specular and diffuse reflection models, is commonly employed. In many practical applications, a fully diffuse boundary condition is utilized to simulate the flow field \cite{bird1994molecular}, and this assumption is also adopted in this work.

For the Maxwell boundary condition, it is assumed that at time $t$ and a boundary point $\bx$, if the wall is at rest, the boundary distribution function $f^b(t, \bx, \bv)$ is determined by the solid wall boundary condition for those velocities $\bv$ that satisfy $\boldsymbol{n}\cdot \bv > 0$, where $\boldsymbol{n}$ is the normal vector of the wall pointing towards the inside of the gas. For other velocities, the boundary distribution function is influenced by the internal distribution function. Assuming the solid wall has the velocity $\boldsymbol{u}^W = (u_1^W,u_2^W,u_3^W)^T$ and temperature $T^W$, the boundary condition is given as 
\begin{equation}
    \label{eq:boun_con}
    f^b(t,\bx,\bv) = \begin{cases}
    \begin{aligned}
        &f^W, &\qquad &\boldsymbol{c}^W \cdot \boldsymbol{n}>0,\\
        &f_N(t,\bx,\bv), &\qquad &\boldsymbol{c}^W \cdot \boldsymbol{n}<0,
    \end{aligned}
    \end{cases}
\end{equation}
where $\boldsymbol{c}^W = \bv -\boldsymbol{u}^W$, $f_N(t,\bx,\bv)$ is the internal distribution function and $f^W$ represents the Maxwellian distribution function determined by the wall velocity $\boldsymbol{u}^W$ and temperature $T^W$ as 
\begin{equation}
    \label{eq:fw}
    f^W\left(\rho^W,\boldsymbol{u}^W,T^W\right) = \frac{\rho^W}{\sqrt{2\pi T^W}^3}\exp\left(-\frac{|\bv-\boldsymbol{u}^W|^2}{2T^W}\right).
\end{equation}
Here, the density $\rho^W$ is decided by the constraint that the normal flux on the wall is zero, which can be expressed as
\begin{equation}
    \label{eq:flux}
    \int_{\bbR^3}\left(\boldsymbol{c}^W  \cdot \boldsymbol{n}\right) f(t,\bx,\bv) \dd\bv = 0.
\end{equation}
By substituting the boundary condition \eqref{eq:boun_con} and Maxwellian distribution function \eqref{eq:fw} into \eqref{eq:flux}. The density $\rho^W$ is calculated as 
\begin{equation}
\label{eq:bound_rho}
    \displaystyle \rho^W = -\frac{\int_{\boldsymbol{c}^W  \cdot \boldsymbol{n}<0}(\boldsymbol{c}^W  \cdot \boldsymbol{n}) f_N(t,\bx,\bv) \dd\bv}{\int_{\boldsymbol{c}^W  \cdot \boldsymbol{n}>0}(\boldsymbol{c}^W  \cdot \boldsymbol{n}) f^W\left(1,\boldsymbol{u}^W,T^W\right) \dd\bv}.
\end{equation}

In the simulation, the Maxwell boundary \eqref{eq:boun_con} with \eqref{eq:bound_rho} is utilized for the microscopic flow problems, and the loss function is specifically designed for this boundary condition which will be introduced in Sec. \ref{sec:loss}.  
 
\subsection{Dimension reduction with flexible auxiliary distribution function}
\label{sec:dim_rec}
In certain practical numerical experiments,  we will encounter scenarios where the dimension of the spatial space $D_x$ is lower than that of the microscopic velocity $D_v=3$. In this case, the dimension reduction model is always utilized to reduce the computational cost. Following the method described in \cite{chu1965kinetictheoretic, yang1995rarefied}, the auxiliary distribution functions are introduced to reduce the dimension of the distribution function for the BGK-type collision model. However, in some applications, such as the Couette flow problem, where there exists macroscopic velocity in a different direction of the spatial space, the dimension of the auxiliary distribution function in the microscopic velocity is larger than that of the spatial space $D_x$.  

Therefore, we enhance this dimension reduction method by incorporating flexible auxiliary distribution functions in this section, in which case, the dimension of the auxiliary distribution function in the microscopic velocity space is reduced to the same as that of the spatial dimension. Precisely, we consider the case $D_x<D_v$. Assuming that the distribution function $f(t, \bx, \bv)$ is uniformly consistent in the dimension $D_x + 1$ to $D_v$ of the spatial space as 
\begin{equation}
\label{eq:rec_f}
    \begin{aligned}
        &f(t,\bx,\bv) = f(t,\tbx,\bv),\qquad \tbx = (x_{1},\cdots,x_{D_x})^T \in \bbR^{D_x},\qquad 
         \frac{\partial f(t,\bx,\bv)}{\partial x_{d}} = 0,\qquad d\in \mathcal{D}, 
    \end{aligned}
\end{equation}
and is plane-symmetric in the dimension $D_{v^{\ast}}$ to $D_v$ of the microscopic velocity space as 
\begin{equation}
\label{eq:plane-sym}
\begin{aligned}
         & f(t,\bx, v_{1}, v_{2}, v_{3}) = f(t,\bx,v_{1}, v_{2}, -v_{3}) , \qquad \text{if}~D_{v^{\ast}}  = 2, \\
         & f(t,\bx, v_{1}, v_{2}, v_{3}) = f(t,\bx,v_{1}, -v_{2}, -v_{3}) , \qquad \text{if}~ D_{v^{\ast}}  = 1.
\end{aligned}
\end{equation}
With \eqref{eq:plane-sym}, we can deduce that the macroscopic velocity $\bu$ satisfies 
\begin{equation}
    u_k = 0,\qquad k \in \hat{\mathcal{D}},
\end{equation}
with 
\begin{equation}
    \begin{aligned}
    \mathcal{D} = \{D_x+1,\cdots,D_v\},\qquad  \hat{\mathcal{D}} = \{D_{v^{\ast}}+1,\cdots,D_v\}, \qquad D_x \leqslant D_{v^{\ast}} \leqslant D_v = 3.
\end{aligned}
\end{equation}
Then, the rigid dimension-reduced distribution functions are defined as 
\begin{equation}
\label{eq:rec_f1}
        g(t,\tbx,\tilde{\bv}) = \int_{\bbR^{D_v-D_x}} f(t,\tbx,\bv) \dd \hat{\bv},\qquad 
        h(t,\tbx,\tilde{\bv}) = \int_{\bbR^{D_v-D_x}} \frac{|\hat{\bv}|^2}{2}f(t,\tbx,\bv) \dd \hat{\bv},
\end{equation}
and the flexible dimension-reduced distribution functions are defined as 
\begin{equation}
    \label{eq:rec_f2}
    s_j(t,\tbx,\tilde{\bv}) = \int_{\bbR^{D_v-D_x}} v_{j} f(t,\tbx,\bv) \dd \hat{\bv},\qquad j = D_x + 1, \cdots, D_{v^{\ast}},
\end{equation}
with 
\begin{equation}
\label{eq:rec_v}
        \tilde{\bv} = \left(v_{1},\cdots,v_{D_x}\right)^T \in \bbR^{D_x},\qquad 
        \hat{\bv} = \left(v_{D_x+1},\cdots,v_{D_v}\right)^T \in \bbR^{D_v-D_x}.
\end{equation}
For simplicity, in subsequent uses of $s_j$ throughout this paper, it is implied that $j = D_x + 1, \cdots, D_{v^{\ast}}$ unless otherwise specified. Here, the number of flexible dimension-reduced distribution functions \eqref{eq:rec_f2} is $D_{v^{\ast}} - D_x$. When $D_{v^{\ast}}=D_x$, flexible dimension-reduced distribution functions are not required, and this new dimension reduction method is the same as in \cite{chu1965kinetictheoretic,yang1995rarefied}. Otherwise, it will further reduce the auxiliary distribution function to the dimension $D_x$ in the microscopic velocity space. 

In this new dimension reduction method, the relationship between the macroscopic variables and the dimension-reduced distribution functions are 
\begin{equation}
    \label{eq:rec_mac_var}
    \begin{aligned}
        &\rho = \int_{\bbR^{D_x}} g(t,\tbx,\tbv) \dd \tilde{\bv},\\
        &\rho u_i = \int_{\bbR^{D_x}} v_{i}g(t,\tbx,\tbv) \dd \tilde{\bv},\qquad i = 1,\cdots,D_x,\\
        &\rho u_j = \int_{\bbR^{D_x}} s_j(t,\tbx,\tbv) \dd \tilde{\bv},\qquad j = D_x+1,\cdots,D_{v^{\ast}},\\
        &E = \frac{3}{2}\rho T + \frac{1}{2}\rho |\bu|^2 = \int_{\bbR^{D_x}}\left( \frac{|\tbv|^2}{2}g(t,\tbx,\tbv)+h(t,\tbx,\tbv)\right) \dd \tilde{\bv}.
    \end{aligned}
\end{equation}
Due to dimension reduction, some higher-order moments cannot be obtained. Therefore, we are limited to calculating only the more significant stress tensor $\sigma_{ij}$ and heat flux $q_i$ 
\begin{equation}
\begin{aligned}
    \sigma_{ij} &= \begin{cases}
        \displaystyle \int_{\bbR^{D_x}} v_{i}v_{j}g(t,\tbx,\tbv) \dd{\tilde{\bv}}-\rho u_iu_j - \delta_{ij}\rho T,\qquad &i,j \leqslant D_x,\\
        \displaystyle \int_{\bbR^{D_x}} v_{i} s_j(t,\bx,\bv)\dd{\tilde{\bv}}-\rho u_iu_j ,&i\leqslant D_x <j \leqslant D_{v^{\ast}},
    \end{cases}\\
    q_{i} &= \displaystyle \int_{\bbR^{D_x}} (v_{i}-u_i)\left(\frac{1}{2}|\tilde{\boldsymbol{c}}|^2g(t,\tbx,\tbv)+h(t,\tbx,\tbv)\right) \dd{\tilde{\bv}} \\ &+\sum_{j=D_x+1}^{D_{v^{\ast}}}\left(u_j^2\int_{\bbR^{D_x}}(v_{i}-u_i) g(t,\tbx,\tbv) \dd{\tilde{\bv}}- 2u_j\int_{\bbR^{D_x}}(v_{i}-u_i) s_j(t,\tbx,\tbv)  \dd{\tilde{\bv}}\right), \qquad i\leqslant D_x,\\
\end{aligned}
\end{equation}
where $\tilde{\boldsymbol{c}} =\tbv - \tilde{\bu} $ with $\tilde{\bu} = (u_1,\cdots,u_{D_x})^T$.

Finally, we derive the governing equation for the reduced system. Multiplying both sides of the Boltzmann equation \eqref{eq:Boltzmann} by $(1,\frac{|\hat{\bv}|^2}{2},v_{D_x+1},\cdots,v_{D_{v^{\ast}}})^T$ and integrate with respect to $\hat{\bv}$ over the domain $\bbR^{D_v-D_x}$,we can obtain the reduced Boltzmann equations:
\begin{equation}
    \label{eq:rec_bgk}
    \begin{aligned}
        &\frac{\partial g}{\partial t} + \tbv \cdot \nabla_{\tbx} g = \frac{1}{\Kn}(g^M-g),\\
        &\frac{\partial h}{\partial t} + \tbv \cdot \nabla_{\tbx} h = \frac{1}{\Kn}(h^M-h),\\
        &\frac{\partial s_j}{\partial t} + \tbv \cdot \nabla_{\tbx} s_j = \frac{1}{\Kn}(s_j^M-s_j),\\
    \end{aligned}
\end{equation}
where $g^M,h^M$ and $s_j^M$ are the related Maxwellian defined as
\begin{equation}
    \label{eq:rec_maxwell}
    \begin{aligned}
        &g^M = \int_{\bbR^{D_v-D_x}} \mathcal{M}[f] \dd \hat{\bv},\qquad 
        &h^M = \int_{\bbR^{D_v-D_x}} \frac{|\hat{\bv}|^2}{2}\mathcal{M}[f] \dd \hat{\bv},\qquad 
        &s^M_j = \int_{\bbR^{D_v-D_x}} v_{j}\mathcal{M}[f] \dd \hat{\bv}.
    \end{aligned}
\end{equation}
The above distribution functions can be explicitly expressed by the macroscopic variables \eqref{eq:rec_mac_var} and we illustrate this with an example where $D_x = 1, D_{v^{\ast}}=2$ as 
\begin{equation}
    \begin{aligned}
        &g^M = \frac{\rho}{\sqrt{2\pi T}}\exp\left(-\frac{|v_{1}-u_1|^2}{2T}\right),\qquad 
        h^M = \frac{(u_2)^2+2T}{2}g^M,\qquad  s_2^M = u_2g^M.
    \end{aligned}
\end{equation}
Then the Maxwell boundary condition introduced in Sec. \ref{sec:boun} can be provided for the reduced system accordingly. With the reduced distribution functions \eqref{eq:rec_f1} and \eqref{eq:rec_f2}, the boundary condition \eqref{eq:boun_con} is changed into
\begin{equation}
    \label{eq:rec_boun_con}
    F^b(t,\tbx,\tbv) = \begin{cases}
    \begin{aligned}
        &F^W, &\qquad &\tilde{\boldsymbol{c}}^W \cdot \tilde{\boldsymbol{n}}>0,\\
        &F_N(t,\tbx,\tbv), &\qquad &\tilde{\boldsymbol{c}}^W \cdot \tilde{\boldsymbol{n}}<0,
    \end{aligned}
    \end{cases}
\end{equation}
 for $F = g,h,s_j$. Here, $\tilde{\boldsymbol{c}}^W = \tbv - \tilde{\bu}^W$ with $\tilde{\bu}^W = (u_1^W,\cdots,u_{D_x}^W)^T$, $\tilde{\boldsymbol{n}}$ is the dimension-reduced normal vector of the wall pointing towards the inside of the gas. $F_N(t,\tbx,\tbv)$ is the internal distribution function and $F^W = g^W,h^W,s_j^W$ is related to \eqref{eq:fw} defined as
\begin{equation}
    \label{eq:rec_boun_maxwell}
    \begin{aligned}
        &g^W = \int_{\bbR^{D_v-D_x}} f^W \dd \hat{\bv},\qquad 
        &h^W = \int_{\bbR^{D_v-D_x}} \frac{|\hat{\bv}|^2}{2}f^W \dd \hat{\bv},\qquad 
        &s^W_j = \int_{\bbR^{D_v-D_x}} v_{j}f^W \dd \hat{\bv}.
    \end{aligned}
\end{equation}
The above distribution function can be explicitly expressed by $\rho^W,\bu^W, T^W$ with $\rho^W$ is calculated as
\begin{equation}
\label{eq:rec_bound_rho}
    \displaystyle \rho^W = -\frac{\int_{\tilde{\boldsymbol{c}}^W \cdot \tilde{\boldsymbol{n}}<0}(\tilde{\boldsymbol{c}}^W \cdot \tilde{\boldsymbol{n}}) g_N(t,\tbx,\tbv) \dd\tbv}{\int_{\tilde{\boldsymbol{c}}^W \cdot \tilde{\boldsymbol{n}}>0}(\tilde{\boldsymbol{c}}^W \cdot \tilde{\boldsymbol{n}}) g^W\left(1,\boldsymbol{u}^W,T^W\right) \dd\tbv}.
\end{equation}
In the following sections, the reduced system \eqref{eq:rec_bgk} and \eqref{eq:rec_maxwell} with the boundary condition \eqref{eq:rec_boun_con} is utilized as the governing equation to design neural representation, which is introduced in the next section. 
% Note that when $D_x = D_v$, if the definition of dimension-reduced distribution functions is still applied, then we obtain $g = f,h=0,p_j = 0$.  In this case, the original Boltzmann equation BGK is a special case of the dimension-reduced system \eqref{eq:rec_bgk}.
\section{The dimension-reduced neural representation}
\label{sec:net}
In this section, the dimension-reduced neural representation (DRNR) is discussed to solve the Boltzmann-BGK equation. We follow three steps, first design the structure of a deep neural network to parameterize the solution to the BGK equation, second, define the loss function to link the neural network and the PDE, and finally train the neural network by minimizing the loss function using optimization algorithms. However, due to the integral-differential property of the BGK equation, it is first semi-discretized in the microscopic velocity space \cite{li2024solving}, where the discretized velocity method (DVM) is utilized, which we will introduce in Sec. \ref{sec:dvm-based_arch}, and then DRNR is proposed in detail, including the network architecture, and the design of the loss function.

\subsection{DVM-based representation}
\label{sec:dvm-based_arch}
Following \cite{li2024solving}, the dimension-reduced model \eqref{eq:rec_bgk} is first discretized in the microscopic velocity space, where a semi-discrete system is obtained as the governing equation to design the neural network. In this case, the input dimension will be reduced greatly, simplifying the neural network's scale and complexity.

First, we abridge $\tbx$ and $\tbv$ defined in \eqref{eq:rec_f} and \eqref{eq:rec_v} as $\bx$ and $\bv$ for clarity in following sections. Assuming the set of the full discrete points in the microscopic velocity space is 
\begin{equation}
\label{eq:V_point}
    \boldsymbol{V} = [\bv_1,\bv_2,\cdots,\bv_{N_v}] \in \bbR^{D_x \times N_v},
\end{equation}
where $N_v$ is the total number of the discrete points, with $\bv_k = \left( v_{1 k},\cdots,v_{D_x k}   \right)^T, 1\leqslant k \leqslant N_v$. Thus, the dimension-reduced distribution functions \eqref{eq:rec_f1} and \eqref{eq:rec_f2} at each discrete velocity point are 
\begin{gather}
\label{eq:dis_f}
    g_k(t,\bx) = g(t,\bx,\bv_k),\qquad
    h_k(t,\bx) = h(t,\bx,\bv_k),\qquad 
    s_{j k}(t,\bx) = s_j(t,\bx,\bv_k), \\
        \label{eq:dis_f_M}
    g^M_k(t,\bx) = g^M(t,\bx,\bv_k),\qquad
    h^M_k(t,\bx) = h^M(t,\bx,\bv_k),\qquad
    s^M_{j k}(t,\bx) = s^M_j(t,\bx,\bv_k).
\end{gather}
Accordingly, the macroscopic variables \eqref{eq:mac_var} can be expressed through discrete distribution functions as \begin{equation}
    \label{eq:dis_mac_var}
    \begin{gathered}
        \rho = \sum_{k=1}^{N_v}g_k\omega_k,\qquad E = \sum_{k=1}^{N_v} \left(  \frac{|\bv_k|^2}{2}g_k+h_k    \right)\omega_k, \\
        \rho \boldsymbol{u}_i = \sum_{k=1}^{N_v}v_{i k}g_k\omega_k,\qquad i = 1,\cdots,D_x, \qquad
        \rho \boldsymbol{u}_j = \sum_{k=1}^{N_v}s_{j k}\omega_k,\qquad j = D_x+1,\cdots,D_{v^{\ast}},
    \end{gathered}
\end{equation}
where $\omega_k, k=1,2,\cdots,N_v$ are the weight of corresponding discrete microscopic velocity points $\bv_k$.
Finally, the semi-discrete form of the reduced BGK equation \eqref{eq:rec_bgk} has the form below 
\begin{equation}
    \label{eq:dis_rec_bgk}
    \begin{aligned}
        &\frac{\partial \boldsymbol{g}}{\partial t} + \boldsymbol{V} \cdot \nabla_{\bx} \boldsymbol{g} = \frac{1}{\Kn}(\boldsymbol{g}^M-\boldsymbol{g}),\\
        &\frac{\partial \boldsymbol{h}}{\partial t} + \boldsymbol{V} \cdot \nabla_{\bx} \boldsymbol{h} = \frac{1}{\Kn}(\boldsymbol{h}^M-\boldsymbol{h}),\\
        &\frac{\partial \boldsymbol{s}_{j}}{\partial t} + \boldsymbol{V} \cdot \nabla_{\bx} \boldsymbol{s}_j = \frac{1}{\Kn}(\boldsymbol{s}_{j}^M-\boldsymbol{s}_{j}),
    \end{aligned}
\end{equation}
where ${\bm F} = {\bm g}, \boldsymbol{g}^M, {\bm h}, \boldsymbol{h}^M, {\bm s}_j, \boldsymbol{s}^M_j$ is discrete dimension-reduced distribution functions, which are defined as
\begin{equation}
\label{eq:dis_f2}
  {\bm F} = [F_1,F_2,\cdots,F_{N_v}] \in \bbR^{1\times N_v}, \qquad F_k = g_k, g^M_k, h_k, h^M_k, s_{jk},s^M_{jk}, \qquad k = 1, \cdots, N_v.
    % \begin{aligned}
    %     \boldsymbol{g} = [g_1,g_2,\cdots,g_{N_v}] \in \bbR^{1\times N_v},\quad 
    %     \boldsymbol{h} = [h_1,h_2,\cdots,h_{N_v}] \in \bbR^{1\times N_v},\quad 
    %     \boldsymbol{s}_j = [s_{j 1},s_{j 2},\cdots,s_{j N_v}] \in \bbR^{1\times N_v}.
    % \end{aligned}
\end{equation}
The semi-discrete dimension-reduced system \eqref{eq:dis_rec_bgk} then is utilized as the governing equations to design the network and loss function in DRNR, which we will introduce in detail in the following sections. 

\subsection{Network architecture}
\label{sec:Net}
In this section, the architecture of the network is introduced. Here, a fully connected neural network is applied to represent the reduced distribution function. The general structure of a $L$-layer fully connected neural network or the feedforward network contains an input layer, $L-1$ hidden layers, and an output layer. Assuming the input variable is $\hat{\bx} \in \mathbb{R}^{m_{0}\times 1}$ and the final output is $\boldsymbol{y}\in \mathbb{R}^{m_{L}\times 1}$, the form of a $L$-layer fully connected neural network can be defined recursively as 
\begin{equation}
    \begin{aligned}
        &\boldsymbol{f}^{[0]}_{\theta}(\hat{\bx}) = \hat{\bx},\\
        &\boldsymbol{f}^{[l]}_{\theta}(\hat{\bx}) = \sigma^{[l]} \circ (\mathbf{W}^{[l]}\boldsymbol{f}_{\theta}^{[l-1]}(\hat{\bx})+{\bm b}^{[l]}),\qquad 1\leqslant l \leqslant L-1,\\
        &\boldsymbol{y} = \boldsymbol{f}^{\rm NN}_{\theta}(\hat{\bx}) = \boldsymbol{f}^{[L]}_{\theta}(\hat{\bx}) =  \mathbf{W}^{[L]}\boldsymbol{f}^{[L-1]}_{\theta}(\hat{\bx})+ {\bm b}^{[L]},
    \end{aligned}
\end{equation}
where $\mathbf{W}^{[l]}\in \mathbb{R}^{m_{l}\times m_{l-1}}, {\bm b}^{[l]}\in \mathbb{R}^{m_{l}\times 1} $, with $m_l$ representing the number of neurons in the $l$-th $(0\leqslant l \leqslant L)$ layer of the neural network. The variable $\theta$ here, denotes the set of trainable parameters, including the weight matrix $\mathbf{W}^{[l]}$ and bias ${\bm b}^{[l]}$. $\sigma^{[l]}$ is the activation function in the $l$-th layer, which is always chosen as a nonlinear function. In this work, the sine and tanh activation functions
\begin{equation}
    \label{eq:act_fun}
    \sigma(x) = \sin(x), \qquad \sigma(x) = \tanh(x)
\end{equation}
are utilized in the simulations due to their superior representation ability \cite{lecun2002efficient, sitzmann2020implicit}.
The symbol ``$\circ$'' means entry-wise operation, as the activation function $\sigma^{[l]}$ is applied on each entry of 
$(\mathbf{W}^{[l]}\boldsymbol{f}_{\theta}^{[l-1]}(\hat{\bx})+{\bm b}^{[l]})$.

To solve the reduced Boltzmann equation \eqref{eq:dis_rec_bgk} using deep neural networks (DNNs), the reduced discrete distribution functions $\bg(t, \bx), \bh(t, \bx)$ and $\bs_j(t, \bx)$ are represented using DNNs. Precisely, the inputs variable $\hat{\bm x}$ here is the spatial variable $\bx$ and time $t$ as $\hat{\bm x} = (t, {\bm x})$, and the discrete distribution functions are taken as outputs  
\begin{equation}
    \label{eq:rep_f}
        {\bm g}(t, \bx) \approx {\bm g}_{\theta}^{\rm NN}(t, \bx), \qquad 
        {\bm h}(t, \bx) \approx {\bm h}_{\theta}^{\rm NN}(t, \bx), \qquad 
        {\bm s}_j(t, \bx) \approx {\bm s}_{j, \theta}^{\rm NN}(t, \bx),(t, \bx), 
\end{equation}
where ${\bm F}_{\theta}^{\rm NN}(t, \bx), {\bm F} = {\bm g}, {\bm h}, {\bm s}_j$ are the $L$-layer fully connect neural network whose dimension of the output is $N_v$, and its $k$-entry of the output is to approximate $F_k$ as 
\begin{equation}
    \label{eq:rep_f_k}
    F_k(t, \bx) \approx F_{\theta, k}^{\rm NN}(t, \bx), \qquad F_k = g_k, h_k, s_{jk}, \qquad k = 1, \cdots, N_v.
\end{equation}

For now, we have introduced the network architecture to represent the dimension-reduced distribution function, and the strategies as multi-scale input and Maxwellian splitting, are applied to enhance the approximation efficiency.

Due to the multi-scale property of the Boltzmann equation, it is believed that the multi-scale networks or adopting multi-scale inputs in the network can significantly enhance the efficiency of the network \cite{xu2019frequency, huang2021solving}. Thus, a similar multi-scale input strategy as in \cite{li2024solving} is utilized to improve the approximation efficiency. Precisely, it is achieved by simply multiplying the input of the network by a sequence of constraints. For all the neural network in \eqref{eq:rep_f} whose input is $\hat{\bx} = (t, \bx) \in \bbR^{D_x + 1}$, the revised input is changed into 
\begin{equation}
    \label{eq:multi_input}
    \hat{\bx}^{\rm multi} = (c_1\hat{\bx}, c_2\hat{\bx},\cdots, c_K\hat{\bx}) \in \bbR^{K(D_x + 1)},
\end{equation}
where $c_i, i = 1,\cdots,K$ are problem dependent constants.

Another strategy utilized here is the Maxwellian splitting \cite{jin2024asymptotic}, which is based on the idea of Micro-Macor decomposition \cite{jin2010micromacro, gamba2019micro}. Since the value of the distribution function at different velocity points can vary significantly, it can lead to difficulties in training the neural network if the output neurons are treated equivalently. To address this issue, the Maxwellian splitting method is utilized here. Precisely, the discrete dimension-reduced distribution function \eqref{eq:dis_f2} is decomposed as 
\begin{equation}
    \boldsymbol{F}(t,\bx) = \boldsymbol{F}^{\text{eq}}(t,\bx) + C\boldsymbol{F}^{\text{neq}}(t,\bx), \qquad {\bm F} = {\bm g}, {\bm h}, {\bm s}_j,
\end{equation}
where $C$ is an adjustable parameter. $\boldsymbol{F}^{\text{eq}}(t, \bx) \in \bbR^{1\times N_v}$ is Maxwellian distribution function
\begin{equation}
\label{eq:f_eq}
        \boldsymbol{F}^{\text{eq}}(t, \bx) = \Big(F_1^{\text{eq}}(t, \bx), F_2^{\text{eq}}(t, \bx), \cdots,  F_{N_v}^{\text{eq}}(t, \bx)\Big)
\end{equation}
with 
\begin{equation}
    \label{eq:f_k_eq}
    F_k^{\text{eq}}(t, \bx)=  \frac{\rho_F^{\ast}(t,\bx)}{\sqrt{2\pi T_F^{\ast}(t,\bx)}^{D_x}}\exp\left(-\frac{|\bv_k-\boldsymbol{u}_F^{\ast}(t,\bx)|^2}{2T_F^{\ast}(t,\bx)}\right), \qquad k = 1, \cdots, N_v. 
\end{equation}
Here $\rho_F^{\ast}, \boldsymbol{u}_F^{\ast} = (u_{F,1}^{\ast},\cdots,u_{F, D_x}^{\ast})^T, T_F^{\ast}$ are the pseudo-macroscopic variables used as intermediate variables to calculate the equilibrium distribution function, and they do not correspond to the actual macroscopic variables associated with $F(t,\bx)$. 

\begin{figure}[!htpb]
    \centering
    \subfigure[Sketch of DRNR]{
                \includegraphics[width = 0.45\textwidth]{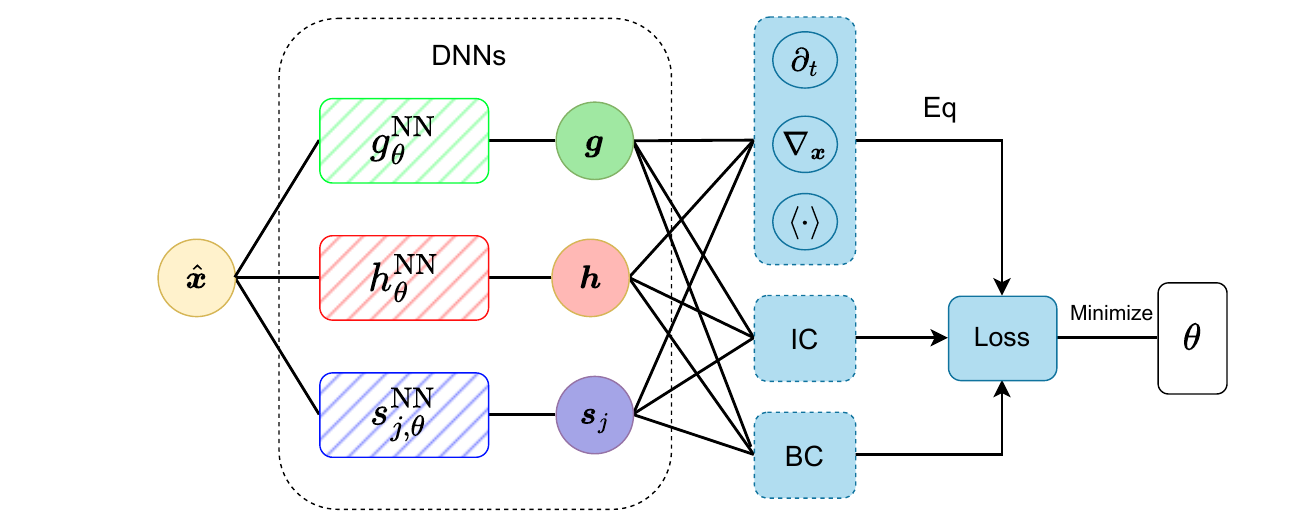}
    }\hfill
    \subfigure[Network architecture]{               
                \includegraphics[width = 0.45\textwidth]{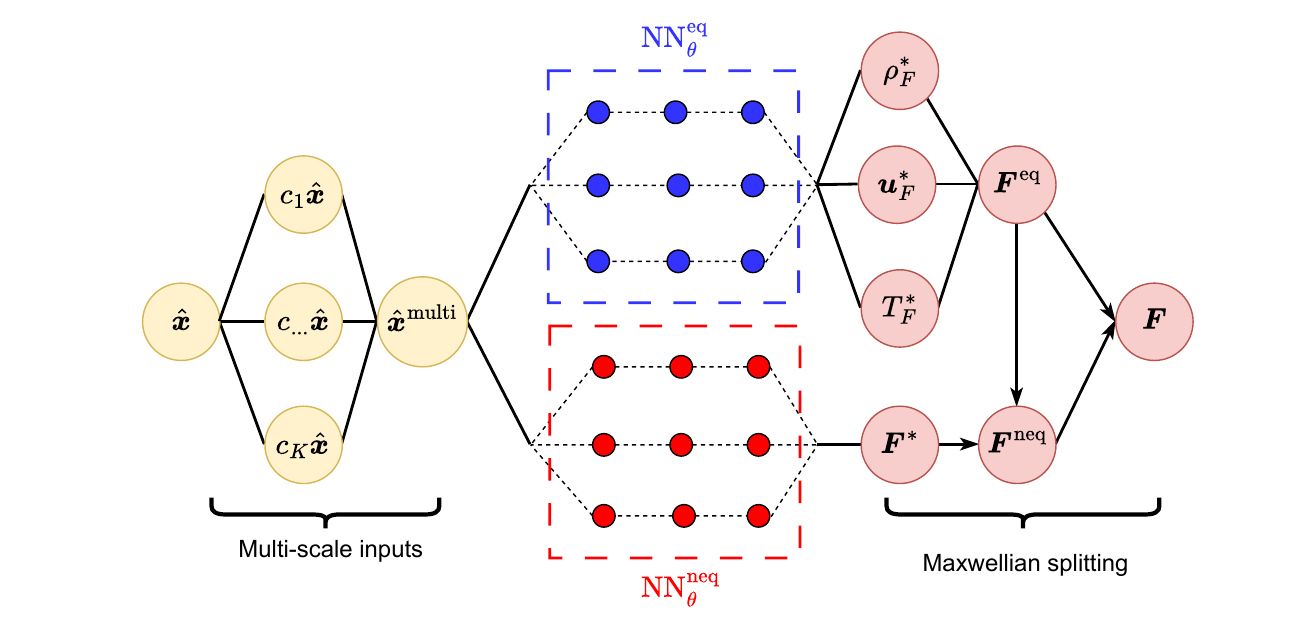}
                \label{fig:net_arch}
    } 
        \caption{Schematic of DRNR for solving the Boltzmann-BGK equation. (a) Sketch of DRNR.        
        In the neural network, the inputs are the spatial space $\bx$ and time $t$. The shaded regions indicate the architecture of the neural networks $g^{\text{NN}}_{\theta},h^{\text{NN}}_{\theta},s^{\text{NN}}_{j,\theta}$, which are used to approximate the dimension-reduced distribution functions ${\bm g}, {\bm h}$, and ${\bm s}_j$. The detailed structure of the neural networks in the shaded region is shown in Fig. \ref{fig:net_arch}.      
        The loss function contains three parts the initial condition (IC), boundary condition (BC), and residual of PDE (Eq).  (b) Network architecture. For the discrete dimension-reduced distribution ${\bm g}, {\bm h}$, and ${\bm s}_j$, each of them corresponds to two separate neural networks, and has its pseudo macro variables. The multi-scale inputs and the Maxwellian splitting are utilized to improve the approximation efficiency. }
    \label{fig:network}
\end{figure}

For the non-equilibrium distribution $\boldsymbol{F}^{\text{neq}}(t, \bx)$, it is assumed that it has the following form 
\begin{equation}
\label{eq:f_neq}
    \boldsymbol{F}^{\text{neq}}(t,\bx) =(F^{\text{neq}}_1(t, \bx),\cdots,F^{\text{neq}}_{N_v}(t, \bx)), \qquad F^{\text{neq}}_k(t, \bx) = F^{\text{eq}}_k(t,\bx) F^{\ast}_k(t,\bx),
\end{equation}
where 
\begin{equation}
    \label{eq:f_neq_star}
    {\bm F}^{\ast}(t,\bx) = (F^{\ast}_1(t, \bx),\cdots,F^{\ast}_{N_v}(t, \bx)) \in \bbR^{1\times N_v}.
\end{equation}
Experimental experience indicates that this can ensure similar magnitudes of ${\bm F}^{\rm eq}(t, \bx)$ and ${\bm F}^{\ast}(t, \bx)$, reducing the difficulties in the training process. Consequently, the pseudo macroscopic variables $(\rho_F^{\ast},\boldsymbol{u}_F^{\ast}, T_F^{\ast})$ and the non-equilibrium discrete distribution function ${\bm F}^{\ast}$ are adopted as the output variables, with two separate neural networks utilized to train them respectively as in \eqref{eq:network_F}.
\begin{equation}
\label{eq:network_F}
\begin{aligned}
    &(\rho_F^{\ast},\boldsymbol{u}_F^{\ast}, T_F^{\ast}) \approx  {\bm F}^{\rm eq, NN}_{\theta}(t, \bx) \triangleq  {\rm NN}_{\theta}^{\text{eq}}(t,\bx) \in \bbR^{D_x+2},\\
    &\boldsymbol{F}^{\ast}(t, \bx) \approx {\bm F}^{\ast, \rm{NN}}_{\theta}(t, \bx) \triangleq  {\rm NN}_{\theta}^{\text{neq}}(t,\bx) \in \bbR^{D_v}, \qquad {\bm F} = {\bm g}, {\bm h}, {\bm s}_j.
\end{aligned}
\end{equation}
 We want to emphasize that for each ${\bm F} = {\bm g}, {\bm h}, {\bm s}_j$, it corresponds to two separate neural networks and has its pseudo-macroscopic variables. Therefore, the total number of the network is  $2(2 + D_{v^{\ast}} - D_x)$. The detailed structure of the neural networks in the shaded region is shown in Fig. \ref{fig:net_arch}.
\begin{remark}
    If we follow the definition in \eqref{eq:rec_f1} and \eqref{eq:rec_f2}, the dimension-reduced distribution functions $g, h$ and $s_j$ are not independent. In this case, a single neural network is required to output the multiple distribution functions, which increases both the output complexity and the optimization difficulty. Therefore, multiple independent neural networks are employed here to represent each dimension-reduced distribution function as shown in \eqref{eq:rep_f}. Independent parameter adjustment allows adaptation to distinct function characteristics, facilitating easier optimization and faster convergence. For each dimension-reduced distribution function, the corresponding equilibrium distribution \eqref{eq:f_eq} is then computed using the same method \eqref{eq:f_k_eq}, where $(\rho_F^{\ast},\boldsymbol{u}_F^{\ast}, T_F^{\ast})$ are the unknown parameters determined by the training process.

    % The equilibrium distribution function \eqref{eq:f_eq} is then computed using the same formula \eqref{eq:f_k_eq} for each independent network as the neural network can accommodate subtle differences during training.
\end{remark}

Once the neural network is built, it is necessary to set a suitable loss function to link the neural network with the PDE. The loss function usually needs to combine the residual of the equation, the initial and boundary conditions, and some additional properties related to PDE, which we will introduce in detail in the next section.

\subsection{Loss function}
\label{sec:loss}
In this section, a specially designed loss function with particular focus on the Maxwell boundary condition is discussed. Generally speaking, the loss function contains the residual of the equation $L_{\text{Eq}}$, the constraints on the 
initial condition $L_{\text{IC}}$ and the boundary condition $L_{\text{BC}}$  as penalty terms 
\begin{equation}
\label{eq:loss}
    L_{\text{loss}} = L_{\text{Eq}}+ \lambda_{1}L_{\text{IC}}+\lambda_{2}L_{\text{BC}},
\end{equation}
where $\lambda_{1}$ and $\lambda_{2}$ are the problem-dependent penalty weights. We will introduce the three parts in detail.

For the time-dependent problem, assuming that the computational domain of the system \eqref{eq:dis_rec_bgk} spans both time and physical space, denoted as $\mathcal{T} \times \Omega$, the training data set is randomly sampled in this domain, where Monte Carlo random sampling is utilized, and uniform random points are generated. For the dimension-reduced system \eqref{eq:dis_rec_bgk}, the residual of the equation $\Leq$ contains three parts as 
\begin{equation}
\label{eq:Loss_eq}
        L_{\text{Eq}} = L^g_{\text{Eq}}+L^h_{\text{Eq}}+\sum_{j=D_x+1}^{D_{v^{\ast}}} L^{s_j}_{\text{Eq}},\qquad  L^F_{\text{Eq}} = \frac{1}{N_{\text{Eq}}}\sum_{i=1}^{N_{\text{Eq}}}\sum_{k=1}^{N_v} \left(R_k^F(t_i,\bx_i)\right)^2, \qquad F = g, h, s_j, 
\end{equation}
where $N_{\rm Eq}$ is the number of the random points, and $R^F_k(t_i,\bx_i)$ is the residual of the semi-discrete reduce BGK equation \eqref{eq:dis_rec_bgk}
\begin{equation}
\label{eq:Loss_residual}
    R^F_k(t_i, \bx_i) = \frac{\partial F^{\text{NN}}_{k,\theta}(t_i, \bx_i)}{\partial t} + \bv_k\cdot \nabla_{\bx} F^{\text{NN}}_{k,\theta}(t_i, \bx_i) - \frac{1}{\Kn}\Big(F^{M}_{k}(t_i, \bx_i)-F^{\text{NN}}_{k,\theta}(t_i, \bx_i)\Big).
\end{equation}
For the steady-state problem, the computational domain of the system is reduced to $\Omega$, and the training data set is randomly sampled on $\Omega$. In this case, \eqref{eq:Loss_residual} is reduced into 
\begin{equation}
\label{eq:Loss_residual-steady}
    R^F_k(\bx_i) =  \bv_k\cdot \nabla_{\bx} F^{\text{NN}}_{k,\theta}(\bx_i) - \frac{1}{\Kn}\Big(F^{M}_{k}(\bx_i)-F^{\text{NN}}_{k,\theta}(\bx_i)\Big).
\end{equation}

For the constraint on the initial condition $\Lic$, the training points are randomly sampled in the region $\{t=0\}\times\Omega$, and the number of points is labeled $N_{\rm IC}$. Similarly, $\Lic$ contains three parts as 
\begin{equation}
    \label{eq:loss_IC}
    L_{\text{IC}} = L^g_{\text{IC}}+L^h_{\text{IC}}+\sum_{j=D_x+1}^{D_{v^{\ast}}} L^{s_j}_{\text{IC}},\qquad 
     L^F_{\text{IC}} =\frac{1}{N_{\text{IC}}}\sum_{i=1}^{N_{\text{IC}}}\sum_{k=1}^{N_v}\left(F^{\text{NN}}_{k,\theta}(0,\bx_i)-F_k^0(\bx_i) \right)^2,\qquad F = g,h, s_j,
\end{equation}
where $F_k^0(\bx_i) = F^0(\bx_0, \bv_k)$ is the given initial condition. For the steady-state problem, the constraint of the initial condition is not included in the loss function. Thus, $\lambda_1$ is set to zero in \eqref{eq:loss}.

For the constraint on the boundary condition $\Lbc$, the training points are randomly sampled in the region $\mathcal{T}\times \partial \Omega$, and the total number of points is labeled $N_{\rm BC}$. $\Lbc$ still contains three parts as 
\begin{equation}
\label{eq:Loss_BC}
L_{\text{BC}} = L^g_{\text{BC}}+L^h_{\text{BC}}+\sum_{j=D_x+1}^{D_{v^{\ast}}} L^{s_j}_{\text{BC}},
\end{equation}
where
\begin{equation}
\label{eq:loss_b}
    L^F_{\text{BC}} =\frac{1}{N_{\text{BC}}}\sum_{i=1}^{N_{\text{BC}}}\sum_{\substack{{k=1}\\{\bv_k\cdot \boldsymbol{n}}>0}}^{N_v}\Big(F^{\text{NN}}_{k,\theta}(t_i,\bx_i)-F^b_k(t_i,\bx_i)\Big)^2,\quad F = g,h,s_j,
\end{equation}
where $F^b_k(t_i, \bx_i) = F^b(t_i, \bx_i, \bv_k)$ is the given boundary condition. In this work, two different boundary conditions are considered, i.e. the Maxwell boundary condition and the inflow boundary condition. For the Maxwell boundary condition, $F^b(t_i, \bx_i, \bv_k)$ is given in \eqref{eq:rec_boun_maxwell} as 
\begin{equation}
    \label{eq:loss_Fb}
    F^b(t_i, \bx_i, \bv_k) = F^W(t_i, \bx_i, \bv_k), \qquad F^W = g^W, h^W, s_j^W,
\end{equation}
where the velocity $\bu^W$ and temperature $T^W$ of the wall are given by the boundary condition while $\rho^W$ is determined based on the condition of the conservation law \eqref{eq:bound_rho}. Here, the internal distribution function is represented by the output of the neural network, and $\rho^W$ is calculated as 
\begin{equation}
    \label{eq:loss_rho}
    \displaystyle \rho^W(t_i,\bx_i) = -\frac{\sum\limits_{\substack{{k=1}\\{\bv_k\cdot \boldsymbol{n}}<0}}^{N_v}(\bv_k  \cdot \boldsymbol{n}) g^{\rm NN}_{k,\theta}(t_i,\bx_i) w_k }{\sum\limits_{\substack{{k=1}\\{\bv_k\cdot \boldsymbol{n}}>0}}^{N_v}(\bv_k  \cdot \boldsymbol{n}) g_k^W(1, u^W(t_i, \bx_i), T^W(t_i, \bx_i))w_k }.
\end{equation}
For the inflow boundary condition, $F^b(t_i, \bx_i, \bv_k)$ is the given boundary condition. 
\begin{remark}
For the steady-state problem, the inputs of neural networks are only the spatial space $\bx$. For most of the steady-state problems, it conserves the total mass, thus we pre-adjust the output of the neural network ${\bm g}_{\theta}^{\rm NN}$ as 
\begin{equation}
    {\bm g}_{\theta}^{\rm NN,steady}(\bx) = \frac{{\bm g}_{\theta}^{\rm NN}(\bx)}{\rho_{\rm{ave}}},
\end{equation}
where
\begin{equation}
    \rho_{\rm{ave}} = \frac{\sum\limits_{i=1}^{N_{\rm Eq}}\langle {\bm g}_{\theta}^{\rm NN}(\bx_i) \rangle }{N_{\rm Eq}},\qquad \langle {\bm g} \rangle = \sum_{k=1}^{N_v} g_k \omega_k. 
\end{equation}
Here, $\bx_i,i=1,\cdots, N_{\rm Eq}$ are the random points on $\Omega$, identical to the spatial training points utilized in the residual of the equation in \eqref{eq:Loss_residual-steady}. This improvement of the network can constrain the conservation of the total mass of ${\bm g}_{\theta}^{\rm NN,steady}(\bx)$ as
\begin{equation}
\sum\limits_{i=1}^{N_{\rm Eq}} \rho_F(\bx_i) \equiv 1, \qquad  F = {\bm g}_{\theta}^{\rm NN,steady}. 
\end{equation}
 This technique can accelerate the training process compared to adding the total mass conservation condition directly to the loss function according to our experience, though there is no theoretical guarantee. A similar technique can also be found in \cite{jin2023asymptotic}. For clarity, we abridge ${\bm g}_{\theta}^{\rm NN, steady}$ as ${\bm g}_{\theta}^{\rm NN}$ afterwards in the steady-state problem.

Moreover, for the steady-state problem, the computational region for the boundary condition is $\partial \Omega$, and the training points are randomly sampled in this region. In this case, \eqref{eq:loss_b} is reduced into 
    \begin{equation}
\label{eq:loss_b_steady}
    L^F_{\rm{BC}} =\frac{1}{N_{\rm{BC}}}\sum_{i=1}^{N_{\rm{BC}}}\sum_{\substack{{k=1}\\{\bv_k\cdot \boldsymbol{n}}>0}}^{N_v}\Big(F^{\rm{NN}}_{k,\theta}(\bx_i)-F^b_k(\bx_i)\Big)^2, \qquad 
        F_k^b(\bx_i) = F^W(\bx_i, \bv_k), \qquad F = g,h,s_j,
\end{equation}
Moreover, for the spatially one-dimension problem such as the Couette flow and Fourier flow problem, the computational region $\partial \Omega$ only has two points. Thus, $N_{\rm BC}$ in \eqref{eq:loss_b} is set as $2$.  
\end{remark}

The loss function \eqref{eq:loss} is a weighted combination of the initial and boundary conditions as well as the PDE residuals. However, for the DVM-based representation, each velocity discrete point has the same weight. In fact, for the calculation of macroscopic quantities, the distribution function is more important when the relative speed is smaller. The lower bound constrained uncertainty weighting method \cite{huang2021solving, li2024solving} is utilized to assign the weight of each microscopic velocity point $\bv_k$ as 
{\small 
\begin{equation}
    \begin{aligned}
        L^F_{\text{Eq}} &= \frac{1}{N_{\text{Eq}}}\sum_{i=1}^{N_{\text{Eq}}}\sum_{k=1}^{N_v} \left(\frac{1}{\Big((w^F_{\text{Eq}})_k + \varepsilon \Big)}R_k^F(t_i,\bx_i)+\log\Big(1+(w^F_{\text{Eq}})_k )\Big)\right)^2,\\
        L^F_{\text{IC}} &=\frac{1}{N_{\text{IC}}}\sum_{i=1}^{N_{\text{IC}}}\sum_{k=1}^{N_v}\left(\frac{1}{\Big((w^F_{\text{IC}})_k + \varepsilon \Big)}\Big(F^{\text{NN}}_{k,\theta}(0,\bx_i)-F_k^0(\bx_i)\Big)+\log\Big(1+(w^F_{\text{IC}})_k\Big) \right)^2,\\
        L^F_{\text{BC}} &=\frac{1}{N_{\text{BC}}}\sum_{i=1}^{N_{\text{BC}}}\sum_{\substack{{k=1}\\{\bv_k\cdot \boldsymbol{n}>0}}}^{N_v}\left(\frac{1}{\Big((w^F_{\text{BC}})_k + \varepsilon\Big)}\Big(F^{\text{NN}}_{k,\theta}(t_i,\bx_i)-F^W_k(t_i,\bx_i)\Big)+\log \Big(1+(w^F_{\text{BC}})_k\Big) \right)^2,
        F = g,h,s_j, 
    \end{aligned}
\end{equation}
}
where $(w^F_{s})_k \geqslant 0,s = \text{Eq, IC, BC},k = 1, \cdots, N_v$ are the adaptive weights for microscopic velocity point $\bv_k$, which are the trainable parameters in the neural network, and $\varepsilon>0$ is a small number to prevent division by zero. 

Right now, we have introduced the description of this dimension-reduced neural representation (DRNR) for the BGK equation, where the fully connected neural network is adopted to approximate the reduced distribution function with the strategies of multi-scale input and Maxwellian splitting utilized to improve the approximation efficiency. The specially designed loss in Sec. \ref{sec:loss} is adopted for the microscopic flow problem to achieve the final results. This method is then applied to solve the BGK equation in Sec. \ref{sec:num}, where its efficiency is validated by the microscopic flow problem.

\section{Numerical examples}
\label{sec:num}
In this section, several numerical experiments are performed to validate the proposed DRNS method for solving the Boltzmann-BGK equation. The numerical experiments include three steady-state problems: 1D Couette flow, Fourier flow, and the 2D rectangular duct flow, as well as one time-dependent problem: the 2D in-out flow. For the steady state problem, the input of the network only requires the spatial variable $\bx$, while for the time-dependent problem, the input of the network requires both the time $t$ and the spatial variable $\bx$.

In the simulation, the Adam optimizer \cite{kingma2014adam} with an initial learning rate of $\eta_0 = 0.001$ is utilized, with the warm restart techniques \cite{loshchilov2017sgdr} incorporated. Precisely, in the $i$-the training step, the learning rate decays according to a cosine annealing schedule as 
\begin{equation}
    \eta_i = \frac{1}{2}\eta_0\left( 1 + \cos \left(\frac{i}{T_{\max}}\right)\right).
\end{equation}
The neural networks comprise five hidden layers, each with 80 neurons and the output dimension is $N_v = \prod\limits_{i=1}^{D_x} N_i$ with $N_i = 100$ and $D_x$ the dimension of the spatial space. The multi-scaling constant in \eqref{eq:multi_input} are chosen as $c_i = 4^{(i-1)},i = 1,2,3$. With the dimension reduction method proposed in Sec. \ref{sec:dim_rec}, the dimension of the microscopic velocity space is the same as the spatial space. Therefore, the computational domain for the microscopic velocity space is chosen as $[-10, 10]^{D_x}$. The detailed descriptions of the sampling points and the network inputs in the training data are provided in each example, respectively.

\subsection{1D Couette flow}
\label{sec:couette}
The Couette flow is a classic problem involving solid wall boundary conditions which is also studied in \cite{hu2020numerical,zhang2023simulation,shakhov1969couette}. The setup consists of two infinite parallel plates at $x = \pm \frac{1}{2}$. The left and right plates have the temperature $T^W=1$ and move in the opposite direction along with the plate with the speed $\boldsymbol{u}^W = (0,\pm u_2^W,0)$ and the density $\rho^W$ decided by \eqref{eq:bound_rho}. Driven by the plate motion, the flow eventually reaches a steady state as time approaches infinity, facilitating studies on the effects of shear in fluid dynamics.

When solving this problem with the dimension-reduced BGK model, one must incorporate the distribution functions $g, h, s_2$ as defined in \eqref{eq:loss_Fb}. Given that the distribution functions are symmetric in $z-$axis in the microscopic velocity space, introducing $s_3$ is unnecessary. The boundary conditions for the Couette flow are specified in  \eqref{eq:rec_boun_maxwell} as follows
\begin{equation} 
\label{eq:couette_ini}
    \begin{aligned}
        &\displaystyle g^W(\rho^W,\boldsymbol{u}^W,T^W) = \frac{\rho^W}{\sqrt{2\pi T^W}}\exp\left(-\frac{|v_1|^2}{2T^W}\right),\\
        &\displaystyle h^W(\rho^W,\boldsymbol{u}^W,T^W) = \frac{|\boldsymbol{u}^W|^2+2T^W}{2}g^W,\\
        &\displaystyle s^W_2(\rho^W,\boldsymbol{u}^W,T^W) = (\boldsymbol{u}^W\cdot \boldsymbol{e}_2) g^W.
    \end{aligned}
\end{equation}
Here, $e_i,i=1,2,3$ represents a unit column vector with only the $i$-th element non-zero.

\begin{figure}[!hptb]
\centering
\subfigure[$\rho (u^W_2 = 0.5)$]{
            \includegraphics[width = 0.45\textwidth]{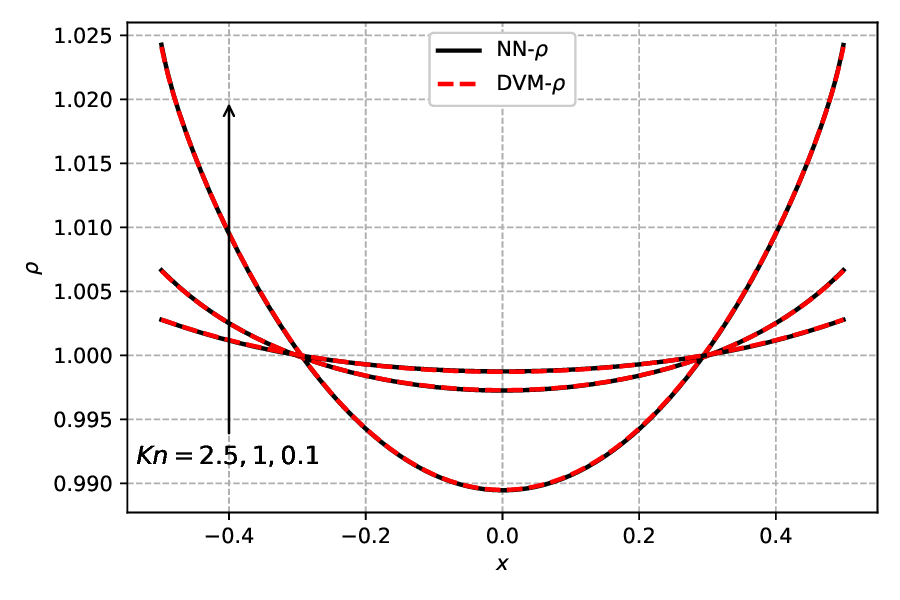}
} \hfill
\subfigure[$u_2 (u^W_2 = 0.5)$]{               
            \includegraphics[width = 0.45\textwidth]{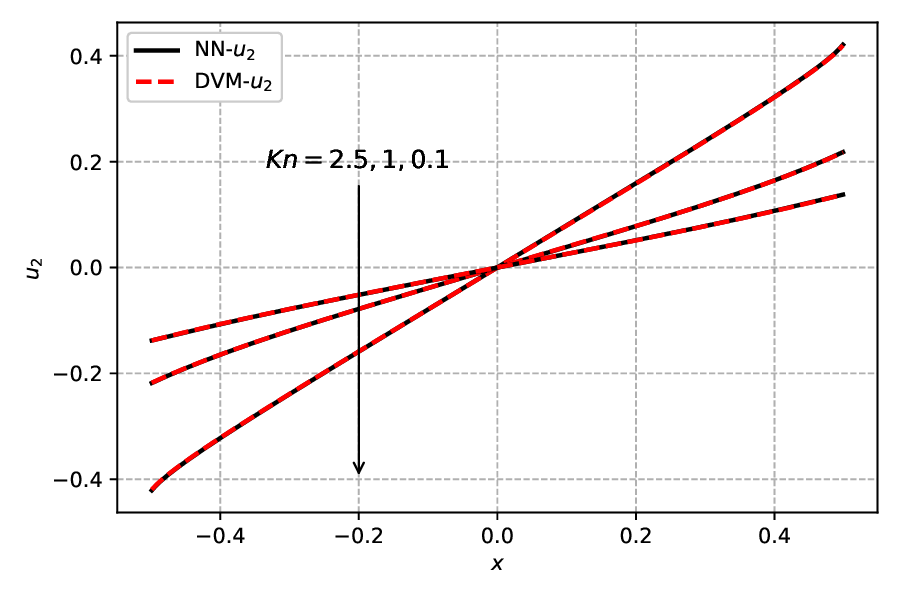}
} \\
\subfigure[$T (u^W_2 = 0.5)$]{
            \includegraphics[width = 0.45\textwidth]{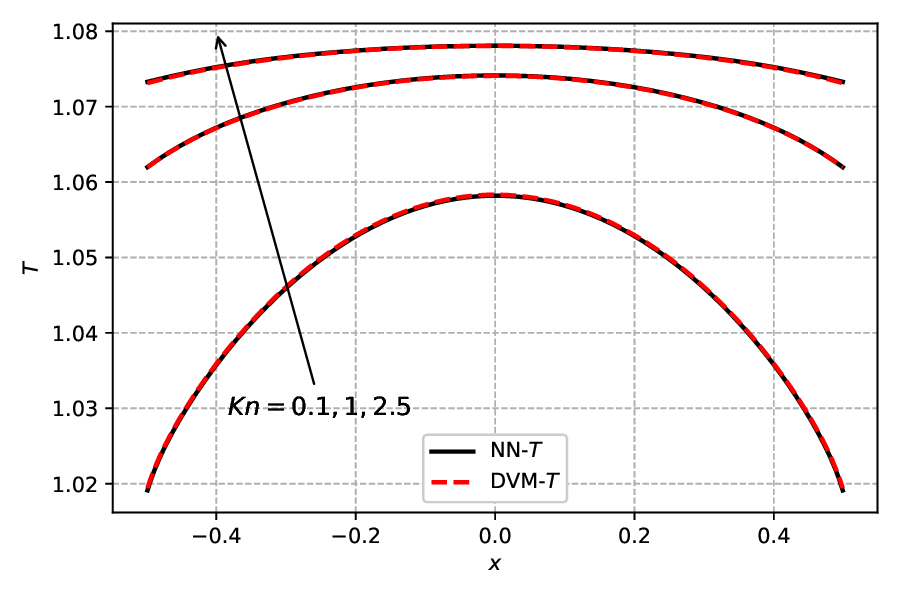}
}\hfill
\subfigure[$q_1 (u^W_2 = 0.5)$]{               
            \includegraphics[width = 0.45\textwidth]{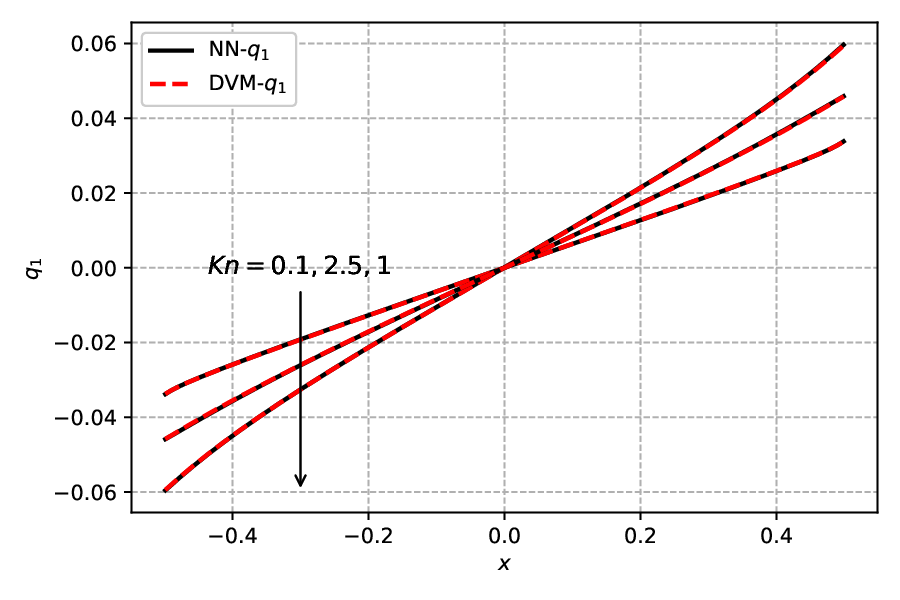}
} 
\caption{(1D Couette flow in Sec \ref{sec:couette}) Numerical solution  of the density $\rho$, macroscopic velocity in $y-$axis $u_2$, the temperature $T$ and the heat flux $q_1$ of the Couette flow at steady state for $\Kn = 0.1, 1$ and $2.5$. Here, the left wall velocity $\bu^W_L = (0,-0.5,0)$ and the right wall velocity $\bu^W_R = (0,0.5,0)$. The black lines are the numerical solution obtained by DRNR, and the dashed red lines represent the reference solution obtained by DVM.}
  \label{fig:couette_cdot5}
\end{figure}

In the simulation, the spatial points are randomly selected with $N_{\rm{PDE}}=500$ in $x\in(-0.5, 0.5)$. For the boundary, one point is chosen fixed at $x  = \pm 0.5$ with $N_{\rm BC} = 2$. The activation function we used is $\sigma(x) = \tanh(x)$ and the total training step is $10,000$. We first set $u^W_2 = 0.5$ and the Knudsen number $\Kn = 0.1, 1$ and $2.5$ are tested. The numerical results of the density $\rho$, the macroscopic velocity in the $y-$axis $u_2$, the temperature $T$, and the heat flux $q_1$ at steady state are plotted in Fig. \ref{fig:couette_cdot5}, with the reference solution obtained by DVM also provided. We find for the density $\rho$, the macroscopic velocity $u_2$, and the temperature $T$, the numerical solution is monotonic with respect to $\Kn$, while the heat flux $q_1$ is increasing with $\Kn$, and then decreasing. This is because when $\Kn = 0$, the Boltzmann equation is reduced into the Euler equation, and the heat flux is zero. When $\Kn$ approaches infinity, in which case there is no collision, the exact solution of the BGK equation can be obtained. We can also calculate the heat flux which is also zero. Therefore, we cannot observe the monotonicity of the heat flux with respect to $\Kn$. Moreover, it shows that for the small Knudsen number $\Kn = 0.1$, the numerical solution matches well with the reference solution, and even when $\Kn$ increases to $\Kn = 1$ and $2.5$, which is far from the equilibrium, the two solutions are also on top with each other. 

\begin{figure}[!hptb]
\centering
\subfigure[$\rho (u^W_2 = 1.0)$]{
            \includegraphics[width = 0.45\textwidth]{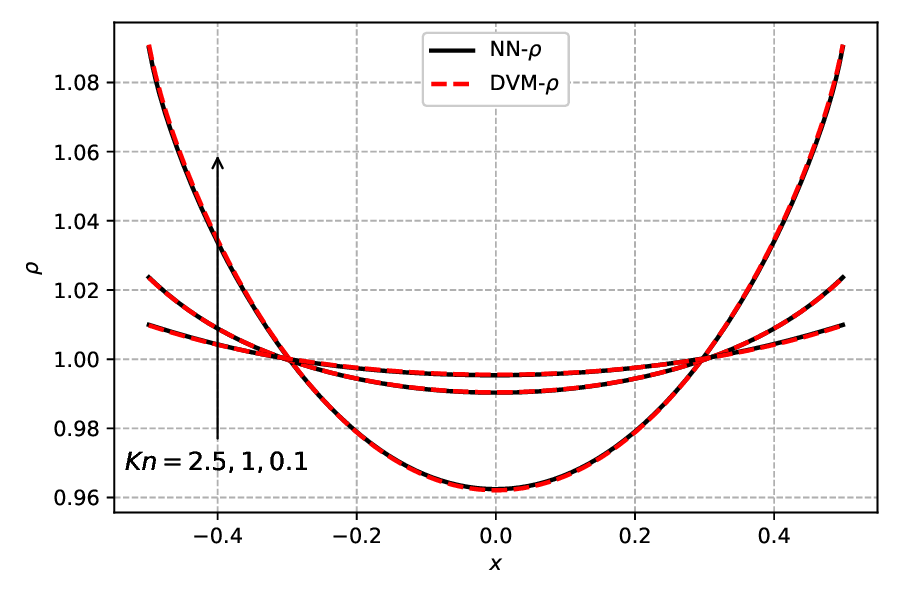}
}\hfill
\subfigure[$u_2 (u^W_2 = 1.0)$]{               
            \includegraphics[width = 0.45\textwidth]{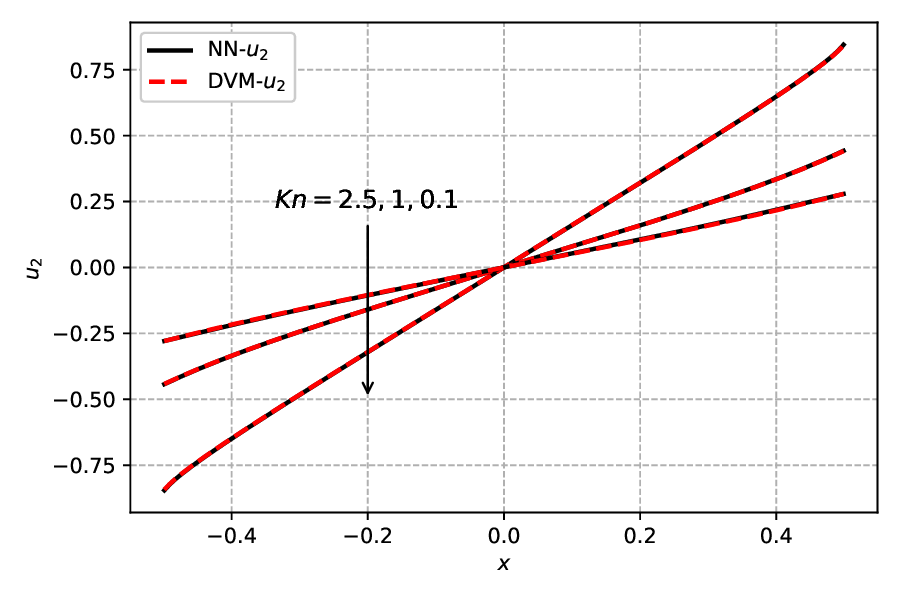}
} \\
\subfigure[$T (u^W_2 = 1.0)$]{
            \includegraphics[width = 0.45\textwidth]{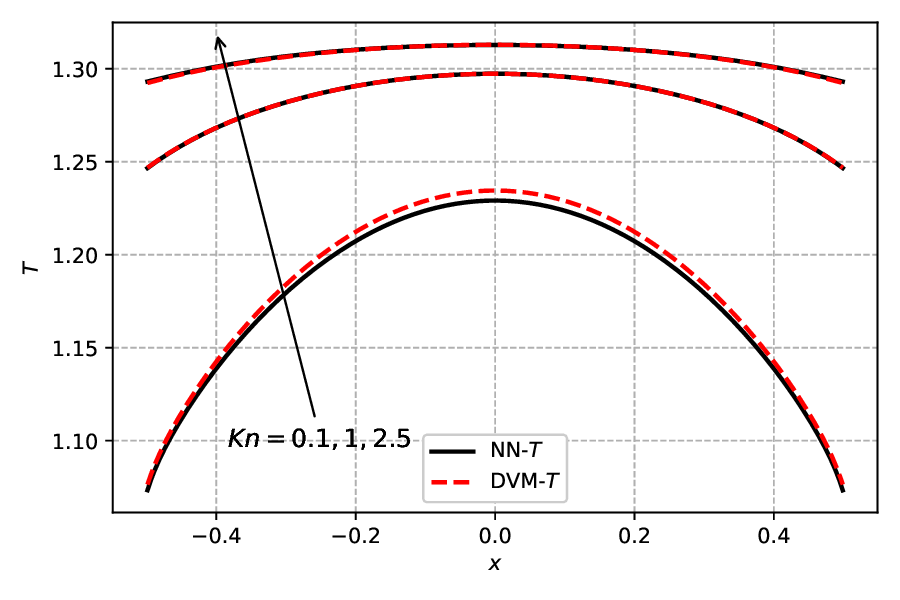}
} \hfill
\subfigure[$q_1 (u^W_2 = 1.0)$]{               
            \includegraphics[width = 0.45\textwidth]{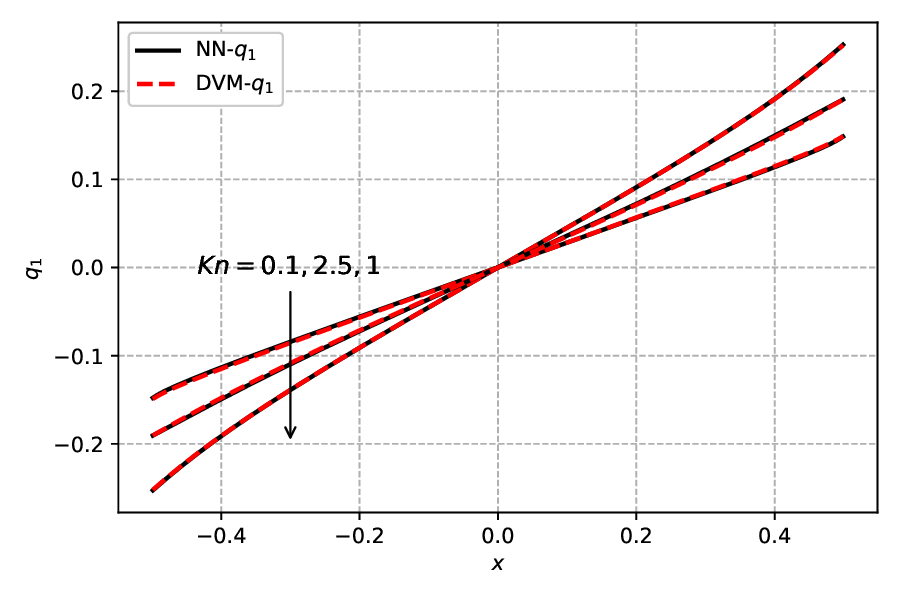}
} 
	\caption{(1D Couette flow in Sec \ref{sec:couette}) Numerical solution  of the density $\rho$, macroscopic velocity in $y-$axis $u_2$, the temperature $T$ and the heat flux $q_1$ of the Couette flow at steady state for $\Kn = 0.1, 1$ and $2.5$. Here, the left wall velocity $\bu^W_L = (0,-1,0)$ and the right wall velocity $\bu^W_R = (0,1,0)$. The black lines are the numerical solution obtained by DRNR, and the dashed red lines represent the reference solution obtained by DVM. }
    \label{fig:couette_c1}
\end{figure}

\begin{figure}[!hptb]
\centering
\subfigure[$\rho (u^W_2 = 2.0)$]{
            \includegraphics[width = 0.45\textwidth]{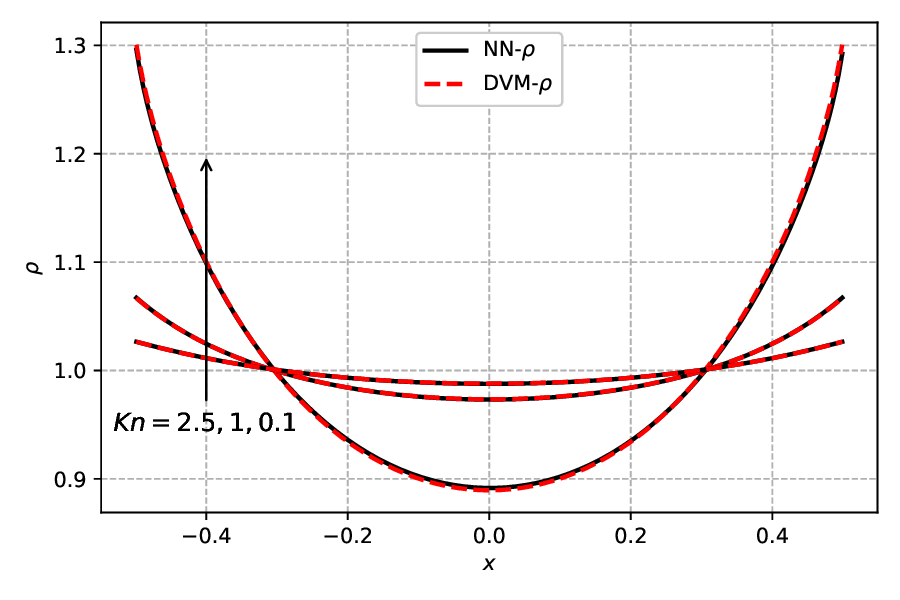}
} \hfill
\subfigure[$u_2 (u^W_2 = 2.0)$]{               
            \includegraphics[width = 0.45\textwidth]{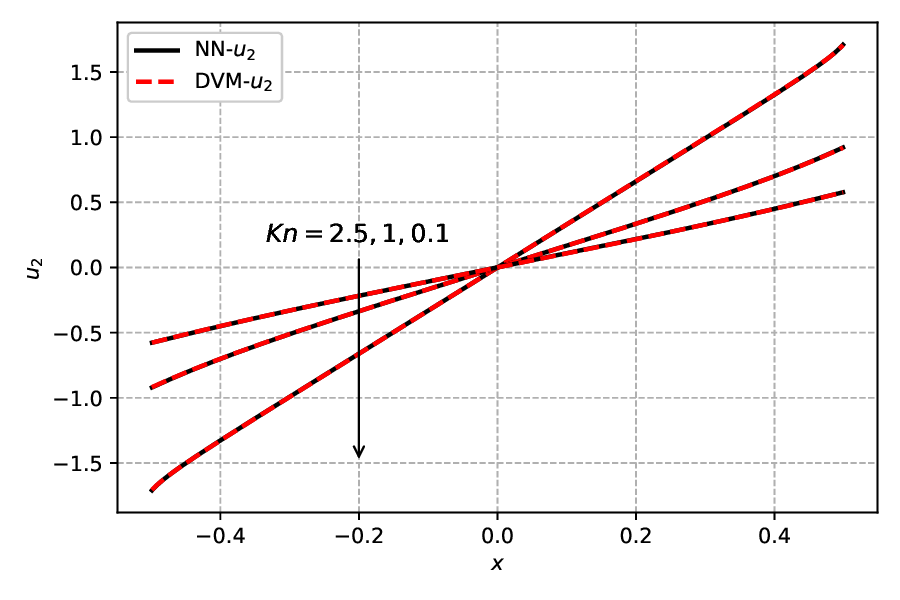}
} \\
\subfigure[$T (u^W_2 = 2.0)$]{
            \includegraphics[width = 0.45\textwidth]{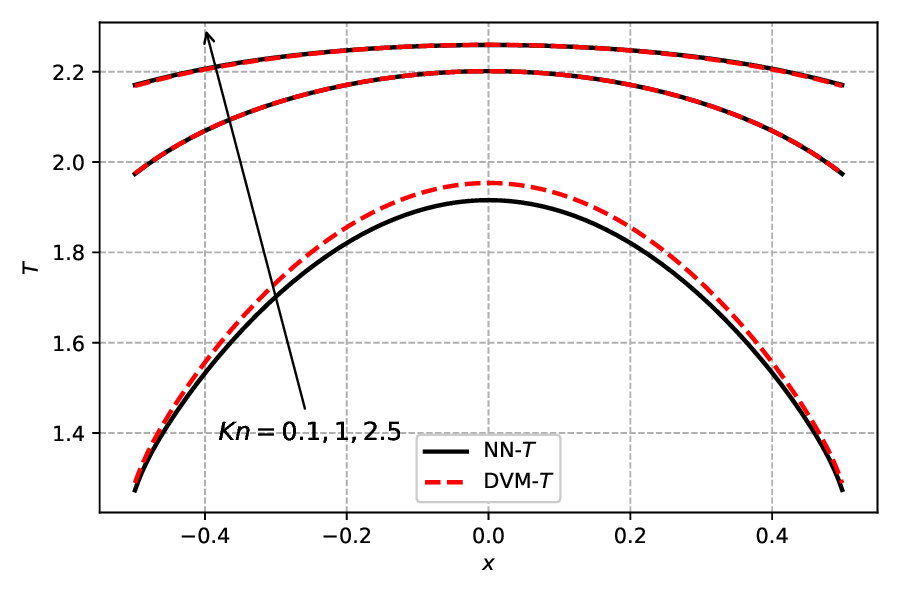}
}\hfill
\subfigure[$q_1 (u^W_2 = 2.0)$]{               
            \includegraphics[width = 0.45\textwidth]{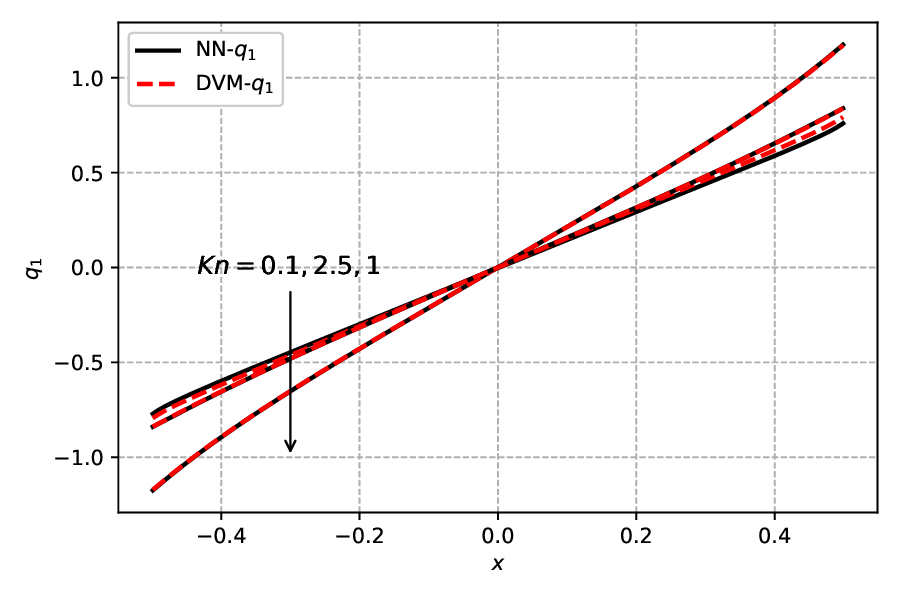}
} 
	\caption{(1D Couette flow in Sec \ref{sec:couette}) Numerical solution  of the density $\rho$, macroscopic velocity in $y-$axis $u_2$, the temperature $T$ and the heat flux $q_1$ of the Couette flow at steady state for $\Kn = 0.1, 1$ and $2.5$. Here, the left wall velocity $\bu^W_L = (0,-2,0)$ and the right wall velocity $\bu^W_R = (0,2,0)$. The black lines are the numerical solution obtained by DRNR, and the dashed red lines represent the reference solution obtained by DVM. }
    \label{fig:couette_c2}
\end{figure}

We increase the wall velocity $u^W_2$ to $1.0$, and the same Knudsen number $\Kn = 0.1, 1$ and $2.5$ are studied. The setting of the network and the training process is the same as that of $u^W_2 = 0.5$. Similarly, the numerical results of the macroscopic variables $\rho$, $u_2$, $T$, and $q_1$ at steady state are illustrated in Fig. \ref{fig:couette_c1}, where the reference solution is also obtained by DVM. The behavior of the numerical solution with respect to $\Kn$ is similar to those of $u^W_2 = 0.1$. Moreover, it shows that for $\rho$, $u_2$ and $q_1$, the numerical solution is well correlated with the reference solution for all three Knudsen numbers. However, for temperature $T$, there is a discrepancy between the numerical solution and the reference solution when $\Kn = 0.1$. But the relative error is less than $1\%$. We deduce that it may be due to the increased nonlinearity of the Boltzmann-BGK equation at smaller Knudsen numbers. For DRNR, it is more difficult to solve the non-linear problem. Finally, we increase $u^W_2$ to $2$, and the numerical results of the same four macroscopic variables at steady state with $\Kn = 0.1, 1$ and $2.5$ are plotted in Fig. \ref{fig:couette_c2}. It is found that, in addition to temperature $T$, there are also some differences between the heat flux $q_1$ when $\Kn = 0.1$. However, the relative error is still quite small. This may be due to the increased $u^W_2$, which strengthens the nonlinearity of the BGK equation.

\begin{figure}[!hptb]
\centering
\subfigure[$\rho$]{
            \includegraphics[width = 0.3\textwidth]{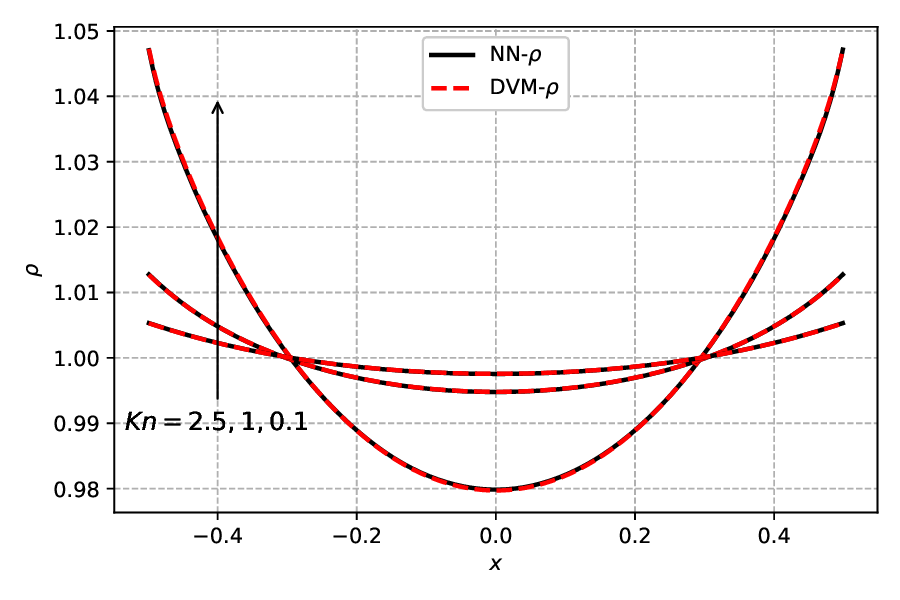}
} \hfill
\subfigure[$u_2$]{
            \includegraphics[width = 0.3\textwidth]{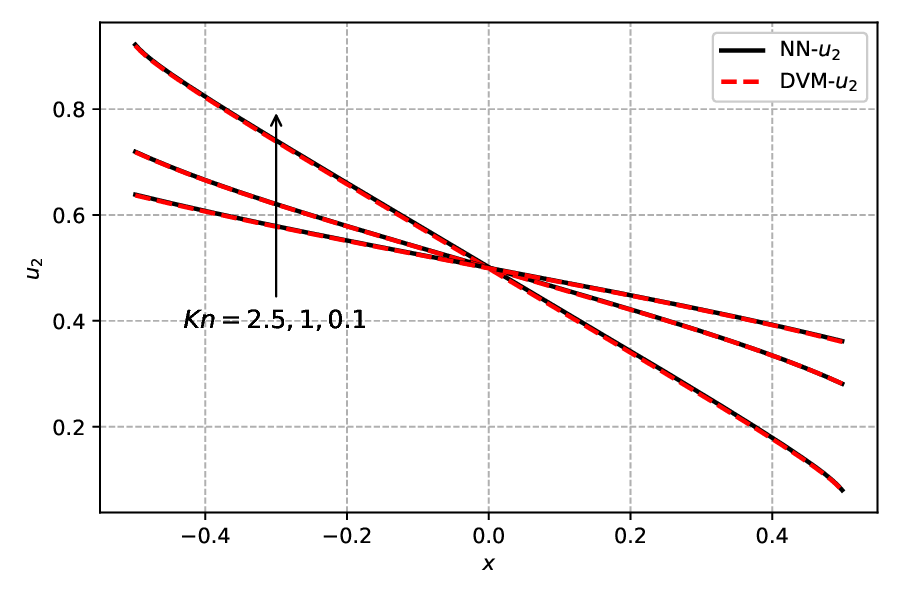}
} \hfill
\subfigure[$u_3$]{               
            \includegraphics[width = 0.3\textwidth]{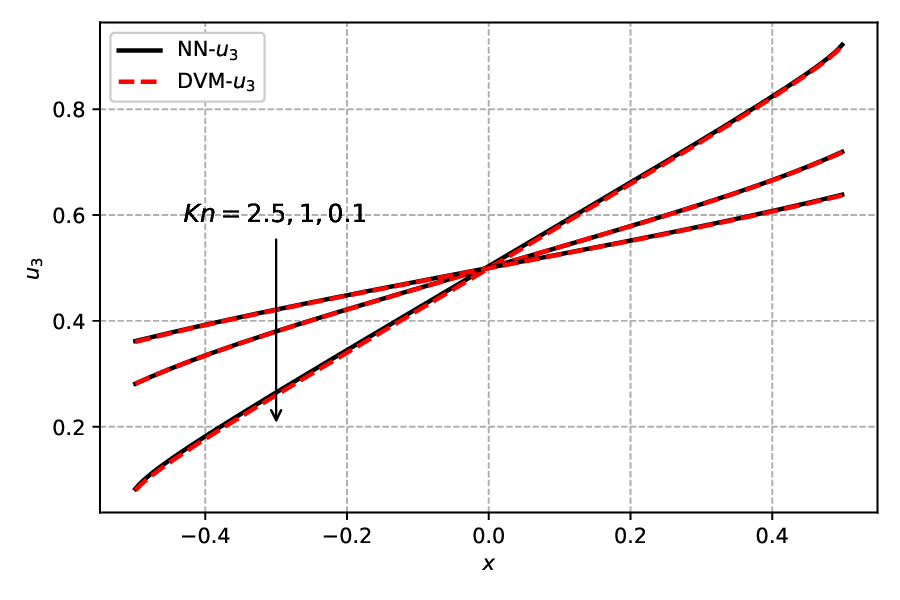}
} \\
\subfigure[$T$]{               
            \includegraphics[width = 0.3\textwidth]{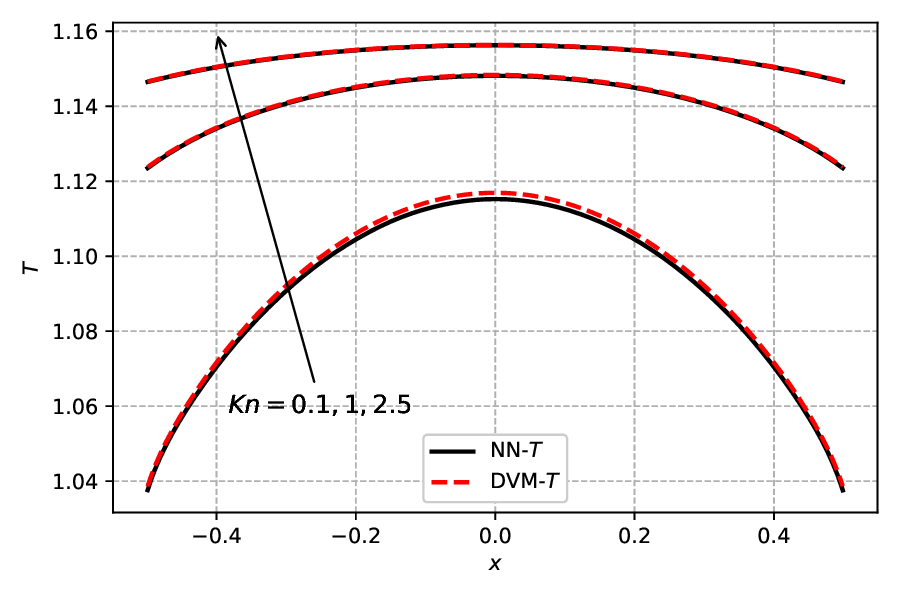}
}\hfill
\subfigure[$\sigma_{13}$]{
            \includegraphics[width = 0.3\textwidth]{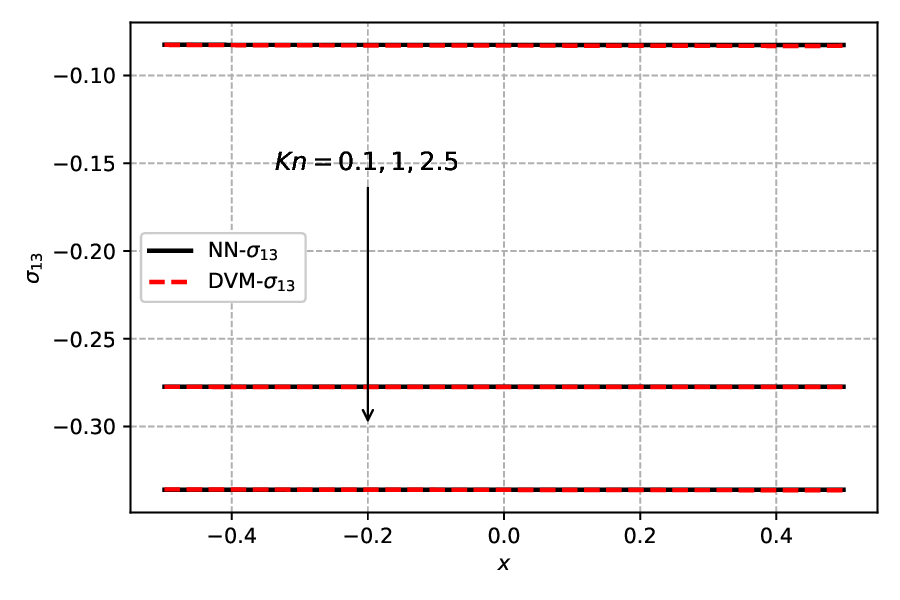}
} \hfill
\subfigure[$q_1$]{               
            \includegraphics[width = 0.3\textwidth]{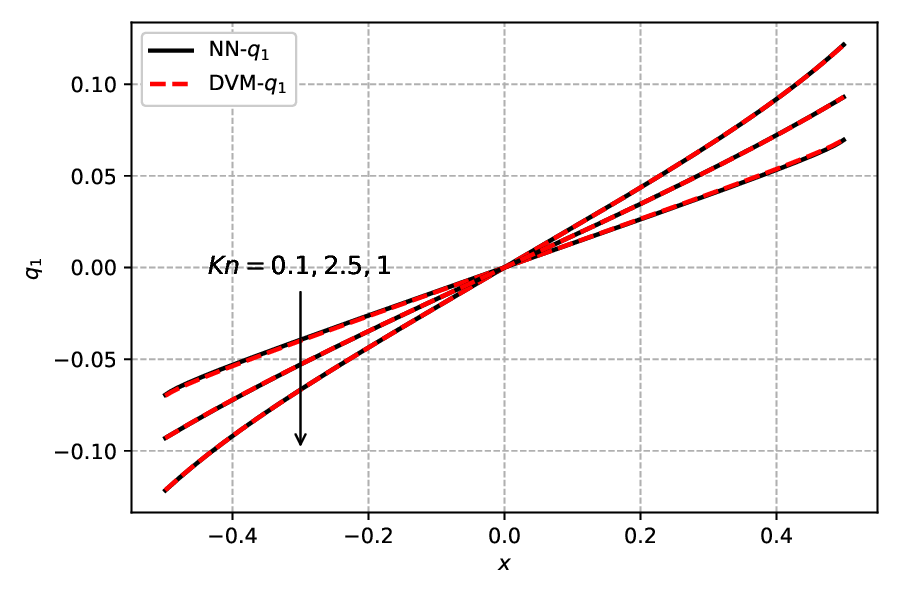}
}
	\caption{(The variant of Couette flow in Sec \ref{sec:couette}) Numerical solution of the density $\rho$, macroscopic velocity $u_2, u_3$, the temperature $T$, the shear stress $\sigma_{13}$ and the heat flux $q_1$ of the variant of Couette flow at steady state for $\Kn = 0.1, 1$ and $2.5$. Here, the left wall velocity $\bu^W_L = (0,1,0)$ and the right wall velocity $\bu^W_R = (0,0,1)$. The black lines are the numerical solution obtained by DRNR, and the dashed red lines represent the reference solution obtained by DVM.}
    \label{fig:C_uz_C1}
\end{figure}

\paragraph{A variant of Couette flow}
Additionally, we explore a variant of Couette flow where the velocities of the left and right plates are in different directions with other parameters the same as that in the last section. The left plate moves at $\bu^W_l = (0,u^W_2,0)$, while the wall velocity of the right plate is $\bu^W_r = (0,0,u^W_3)$. In this case, there is no symmetry in $z-$axis. Consequently, it is necessary to introduce an additional distribution function $s_3$ along with its corresponding boundary condition:
\begin{equation}
    \displaystyle s^W_3(\rho^W,\boldsymbol{u}^W,T^W) = (\boldsymbol{u}^W\cdot \boldsymbol{e}_3) g^W,
\end{equation} 
with $g^W$ defined in \eqref{eq:couette_ini}. 
In the simulation, we set $u^W_2 = u^W_3 = 1$. The network and training parameters are the same as in the last section. The numerical results of the density $\rho$, macroscopic velocity $u_2, u_3$, the temperature $T$, the stress tensor $\sigma_{13}$ and the heat flux $q_1$ at steady state with $\Kn = 0.1, 1$ and $2.5$ are plotted in Fig. \ref{fig:C_uz_C1}, along with the reference solution obtained by DVM. We find that in this variant of Couette flow, the macroscopic velocity $u_3$, and the shear stress $\sigma_{13}$ are non-zero compared to the general Couette flow problem. For all variables, the numerical solution matches well with the reference solution for all three Knudsen numbers. 

\subsection{1D Fourier heat transfer flow}
\label{sec:fourier}
Fourier heat transfer flow is another classical problem, which is studied in \cite{hu2020burnett,zhang2023simulation}. In contrast to Couette flow, two infinite plates are fixed at $x = \pm \frac{1}{2}$, with the plates being stationary but having different temperatures $T^W_l, T^W_r$. The macroscopic velocity of the two plates is $\bu^W = \bz$ with the density $\rho^W$ decided by \eqref{eq:bound_rho}. The flow will also reach a steady state as time goes on. Since the macroscopic velocities of the other two dimensions are zero, the dimension-reduced functions $g$ and $h$ are needed in the Fourier flow problem, and the corresponding boundary conditions are as below 
\begin{equation}
\label{eq:Fourier_ini}
    \begin{aligned}
        &\displaystyle g^W(\rho^W,\boldsymbol{u}^W,T^W) = \frac{\rho^W}{\sqrt{2\pi T^W}}\exp\left(-\frac{|v_1|^2}{2T^W}\right),\\
        &\displaystyle h^W(\rho^W,\boldsymbol{u}^W,T^W) = T^Wg^W, \qquad T^W = T^W_l, T^W_r. 
    \end{aligned}
\end{equation}

\begin{figure}[!hptb]
\centering
\subfigure[$\rho$ with wall temperature \eqref{eq:Fourier_T1}]{
            \includegraphics[width = 0.45\textwidth]{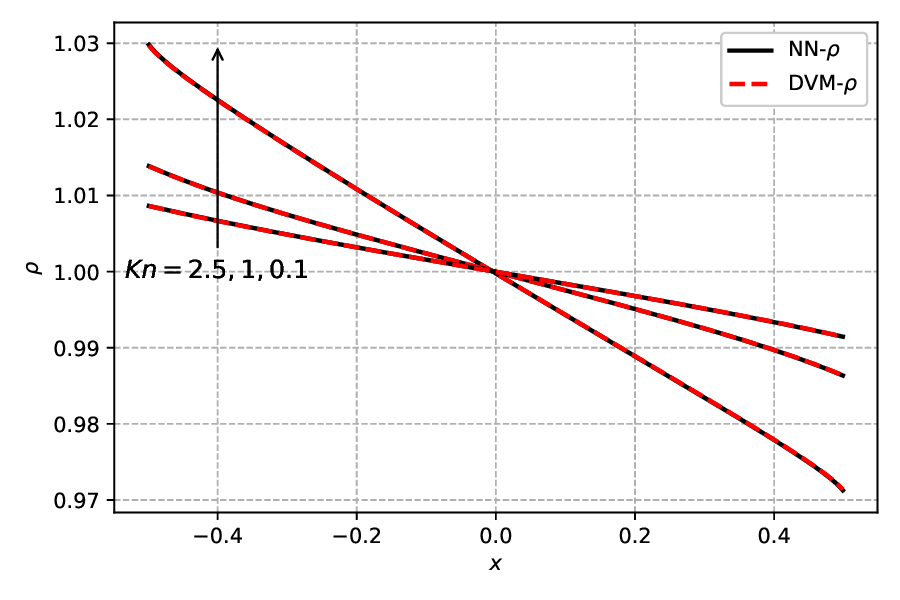}
} \hfill
\subfigure[$T$ with wall temperature \eqref{eq:Fourier_T1}]{               
            \includegraphics[width = 0.45\textwidth]{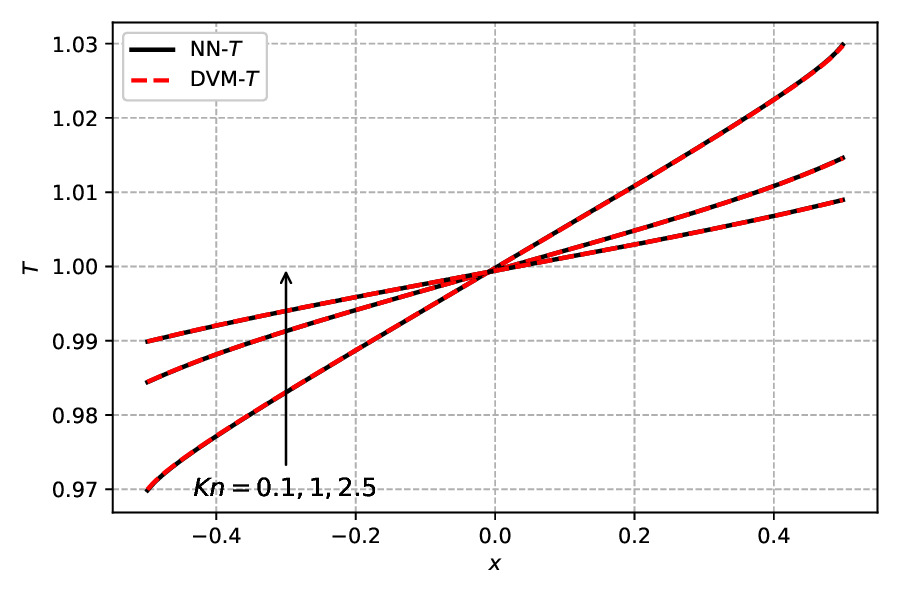}
} \\
\subfigure[$\sigma_{11}$ with wall temperature \eqref{eq:Fourier_T1}]{
            \includegraphics[width = 0.45\textwidth]{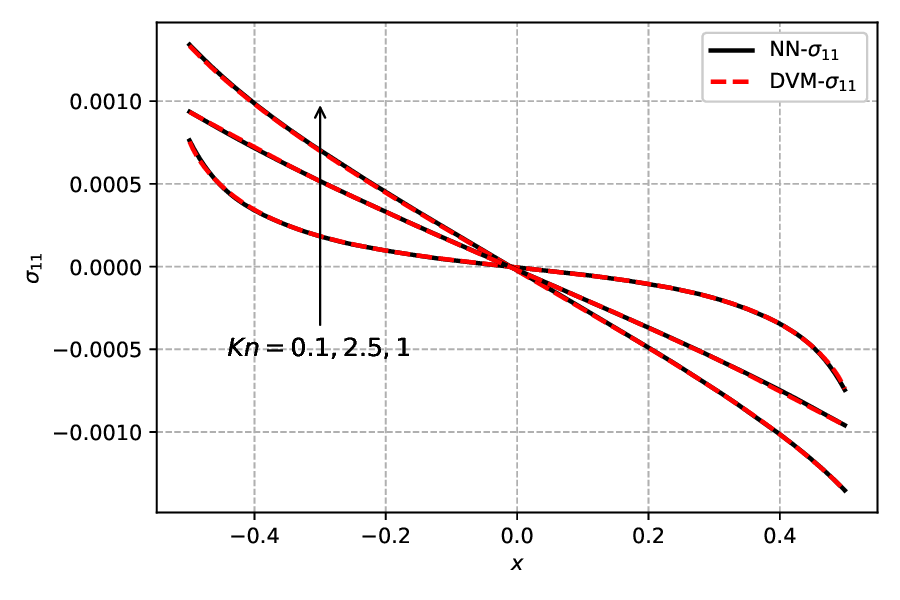}
} \hfill
\subfigure[$q_1$ with wall temperature \eqref{eq:Fourier_T1}]{               
            \includegraphics[width = 0.45\textwidth]{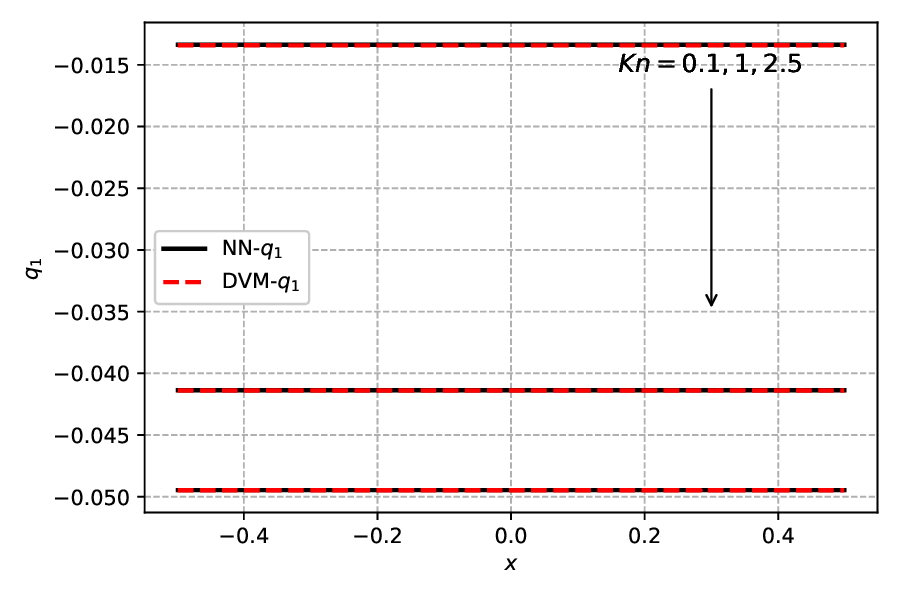}
} 
\caption{(1D Fourier heat transfer flow in Sec \ref{sec:fourier}) Numerical solution  of the density $\rho$, the temperature $T$, the shear street $\sigma_{11}$ and the heat flux $q_1$ of the Fourier flow at steady state for $\Kn = 0.1, 1$ and $2.5$. Here, the left wall temperature is $T_l^W = \frac{T_0 - 10}{T_0}$ and the right wall temperature is $T_r^W = \frac{T_0 + 10}{T_0}$ with $T_0 = 273.15$. The black lines are the numerical solution obtained by DRNR, while the dashed red lines represent the reference solution obtained by DVM.
}
\label{fig:F_t263dot15}
\end{figure}

\begin{figure}[!hptb]
\centering
\subfigure[$\rho$ with wall temperature \eqref{eq:Fourier_T2}]{
            \includegraphics[width = 0.45\textwidth]{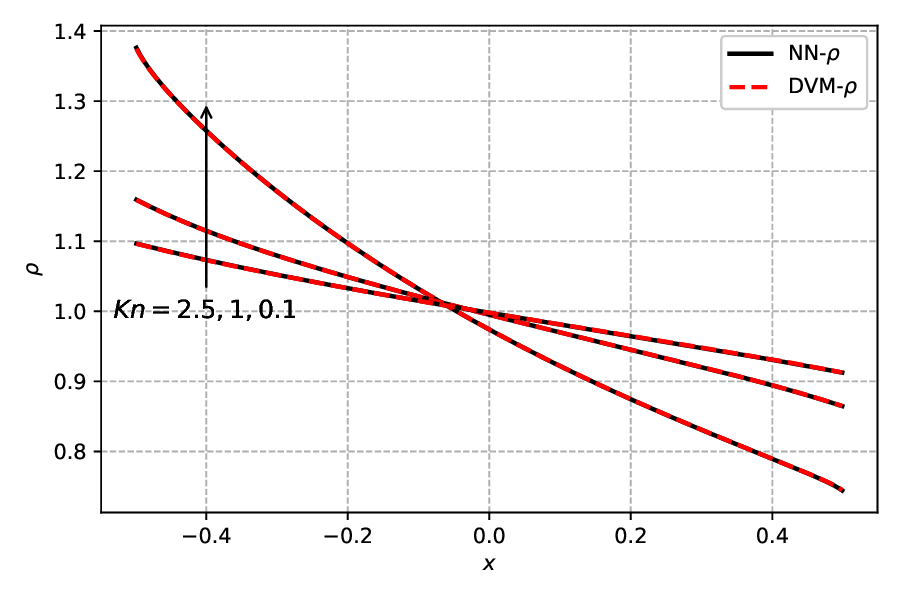}
}
\subfigure[$T$ with wall temperature \eqref{eq:Fourier_T2}]{               
            \includegraphics[width = 0.45\textwidth]{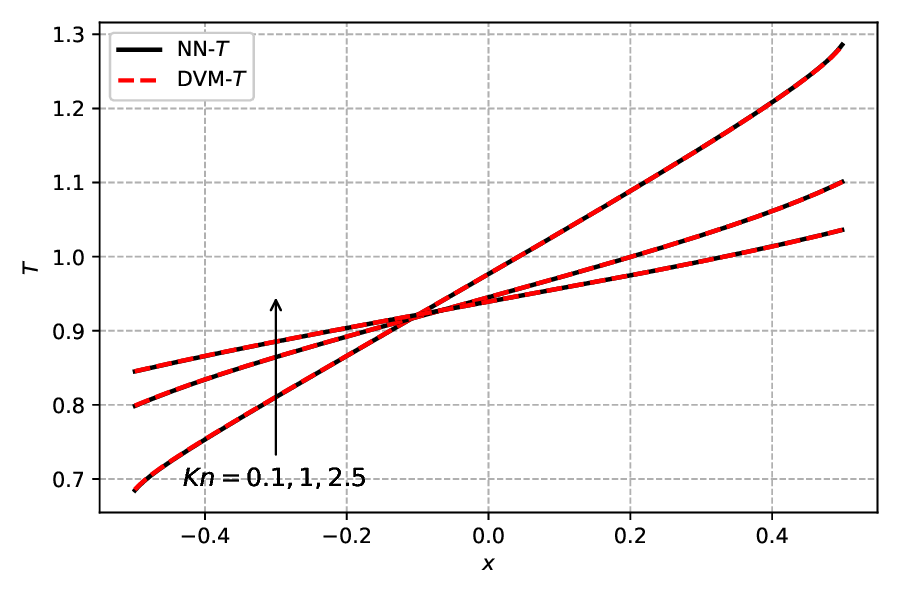}
}       
\subfigure[$\sigma_{11}$ with wall temperature \eqref{eq:Fourier_T2}]{
            \includegraphics[width = 0.45\textwidth]{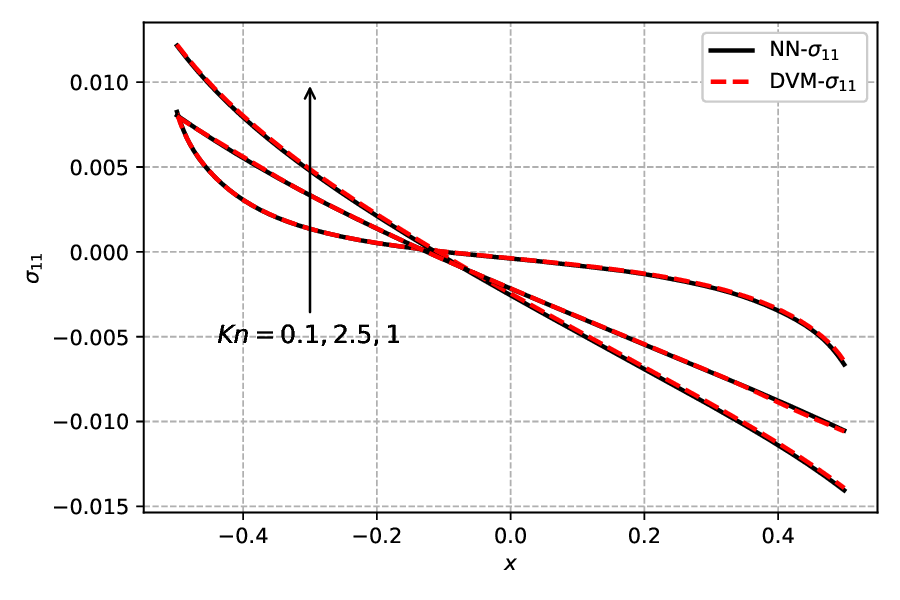}
}
\subfigure[$q_1$ with wall temperature \eqref{eq:Fourier_T2}]{               
            \includegraphics[width = 0.45\textwidth]{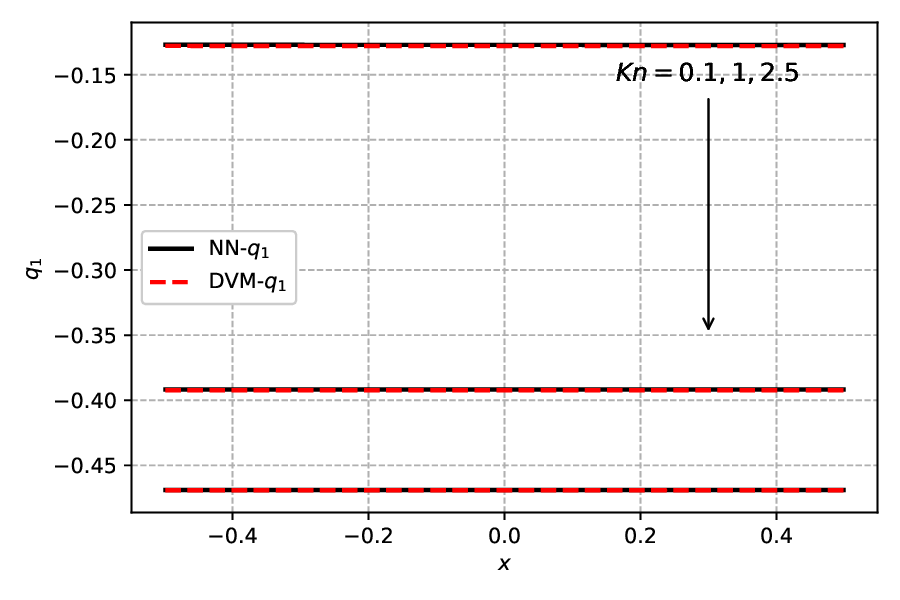}
} 
\caption{
(1D Fourier heat transfer flow in Sec \ref{sec:fourier}) Numerical solution of the density $\rho$, the temperature $T$, the shear street $\sigma_{11}$ and the heat flux $q_1$ of the Fourier flow at steady state for $\Kn = 0.1, 1$ and $2.5$. Here, the left wall temperature is $T_l^W = \frac{T_0 -100}{T_0}$ and the right wall temperature is $T_r^W = \frac{T_0 + 100}{T_0}$ with $T_0 = 273.15$. The black lines are the numerical solution obtained by DRNR, while the dashed red lines represent the reference solution obtained by DVM.}
\label{fig:F_t173dot15}
\end{figure}

We first set the wall temperature as 
\begin{equation}
    \label{eq:Fourier_T1}
    T^W_l = \frac{T_0 - 10}{T_0}, \qquad T^W_r = \frac{T_0 + 10}{T_0}, \qquad T_0 = 273.15.
\end{equation}
For the network, the spatial points are randomly selected with $N_{\rm{PDE}} = 500$ in $x \in (-0.5, 0.5)$. For the boundary, $N_{\rm BC} = 2$ at $x = \pm 0.5$. The activation function we used is $\sigma(x) = \tanh(x)$ and the total training step is $10,000$. The numerical solutions of the density $\rho$, temperature $T$, shear stress $\sigma_{11}$, and the heat flux $q_1$ at steady state for the Knudsen numbers $\Kn = 0.1, 1.0$ and $2.5$ are plotted in Fig. \ref{fig:F_t263dot15}. Here, the reference solution obtained by DVM is also provided. For this small temperature difference, the numerical solution matches well with the reference solution for the three Knudsen numbers. 
For the Fourier flow, the numerical solution of the density $\rho$, the temperature $T$, and the heat flux $q_1$ is monotonic with respect to $\Kn$, while the shear stress $\sigma_{11}$ behaves similarly to that of the heat flux in the Couette flow for the same reason. Then we exchange the difference in the wall temperature as 
\begin{equation}
    \label{eq:Fourier_T2}
    T^W_l = \frac{T_0 - 100}{T_0}, \qquad T^W_r = \frac{T_0 + 100}{T_0}, \qquad T_0 = 273.15.
\end{equation}
With the same network and training process, the numerical results for the four macroscopic variables are shown in Fig. \ref{fig:F_t173dot15} with the reference solution obtained by DVM presented. It indicates that the variations of these macroscopic variables are larger compared to those of \eqref{eq:Fourier_T1}. Even though, the numerical solution and the reference solution are still on top of each other, demonstrating the efficiency of DRNR.

\begin{figure}[!hptb]
\centering
\subfigure[$\rho$ with wall temperature \eqref{eq:Fourier_T3}]{
            \includegraphics[width = 0.45\textwidth]{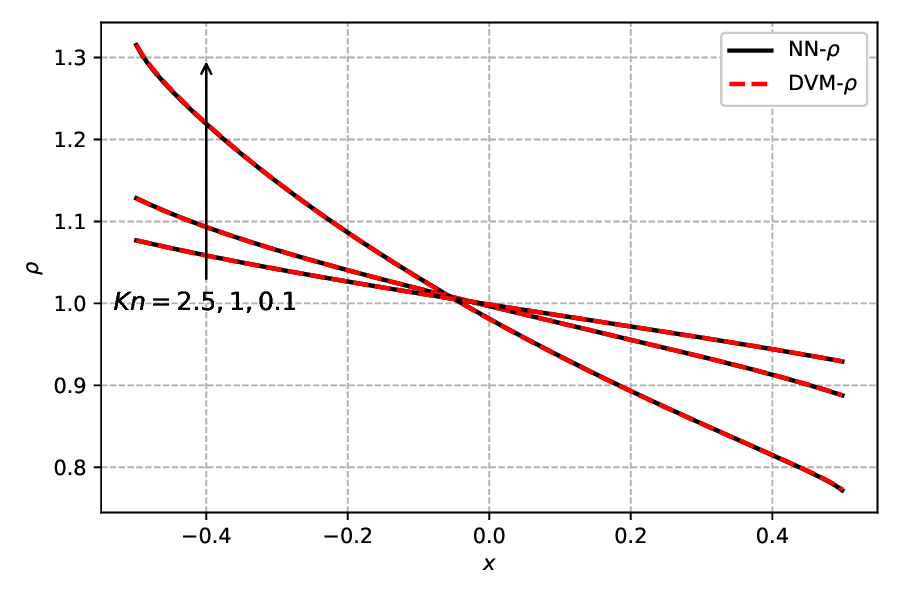}
}
\subfigure[$T$ with wall temperature \eqref{eq:Fourier_T3}]{               
            \includegraphics[width = 0.45\textwidth]{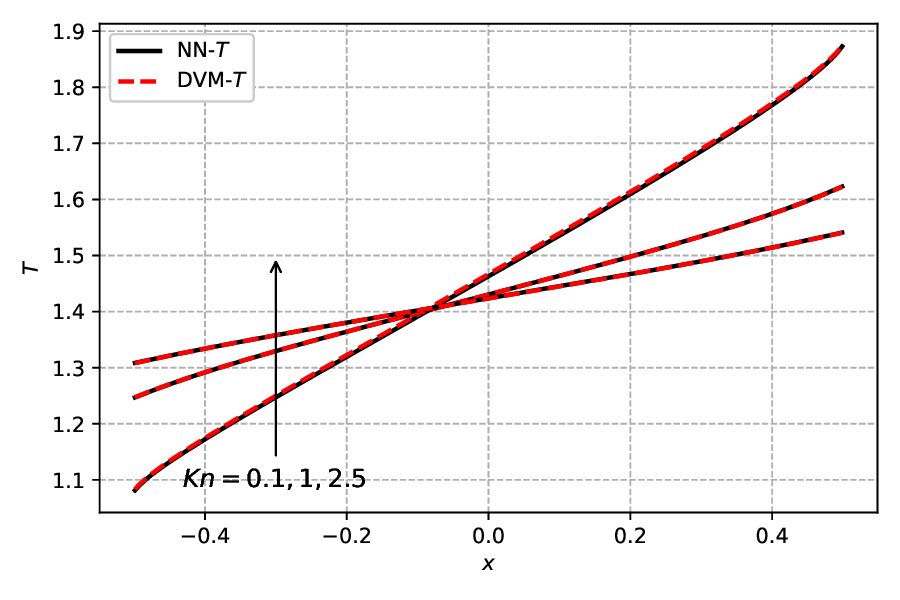}
}       
\subfigure[$\sigma_{11}$ with wall temperature \eqref{eq:Fourier_T3}]{
            \includegraphics[width = 0.45\textwidth]{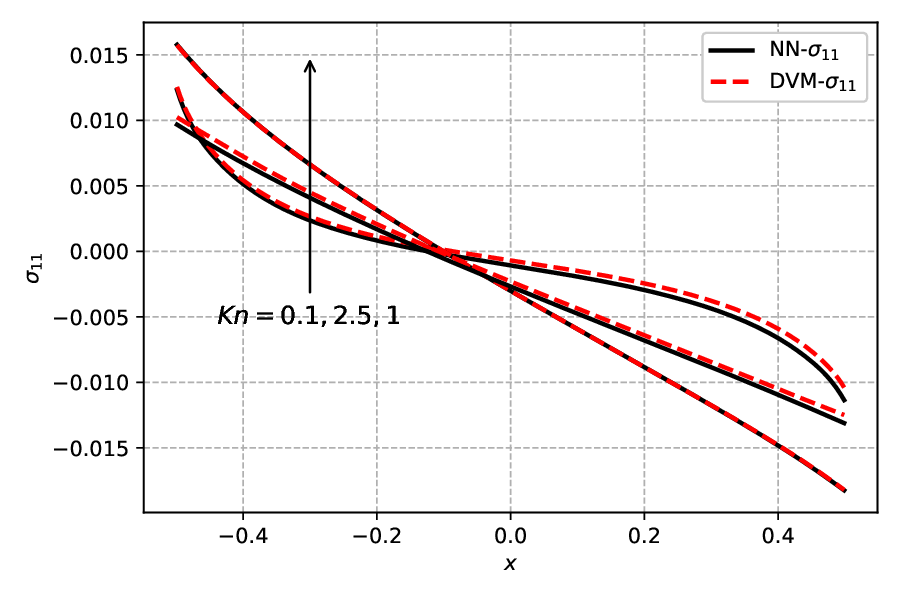}
}
\subfigure[$q_1$ with wall temperature \eqref{eq:Fourier_T3}]{               
            \includegraphics[width = 0.45\textwidth]{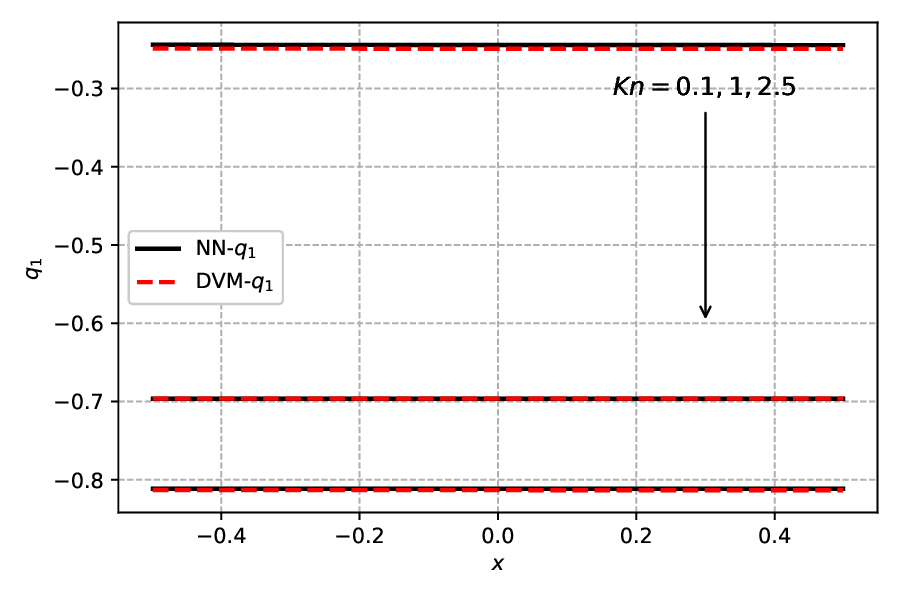}
} 
\caption{
(1D Fourier heat transfer flow in Sec \ref{sec:fourier}) Numerical solution of the density $\rho$, the temperature $T$, the shear street $\sigma_{11}$ and the heat flux $q_1$ of the Fourier flow at steady state for $\Kn = 0.1, 1$ and $2.5$. Here, the left wall temperature is $T_l^W = 1.0$, and the right wall temperature is $T_r^W = 2.0$. The black lines are the numerical solution obtained by DRNR, while the dashed red lines represent the reference solution obtained by DVM.}    
\label{fig:F_t1-2}
\end{figure}

Finally, we make the difference in the wall temperature even larger as
\begin{equation}
    \label{eq:Fourier_T3}
    T^W_l = 1.0, \qquad T^W_r = 2.0.
\end{equation}
The numerical solution and the reference solution obtained by DVM of $\rho$, $T$, $\sigma_{11}$, and $q_1$ at steady state for $\Kn = 0.1, 1.0$ and $2.5$ are illustrated in Fig. \ref{fig:F_t1-2}. We find for $\rho$, $T$, and $q_1$, the numerical solution and the reference solution all match well with each other, while there exists a little difference for $\sigma_{11}$ when $\Kn = 0.1$ and $2.5$. We deduce that it may be due to the non-linearity of the BGK equation when $\Kn = 0.1$, and the non-linearity brought by the large temperature ratio when $\Kn = 2.5$.  

\subsection{2D rectangular duct flow}
\label{sec:cavity}
\begin{figure}[!hptb]
\centering
\includegraphics[width = 0.3\textwidth, clip]{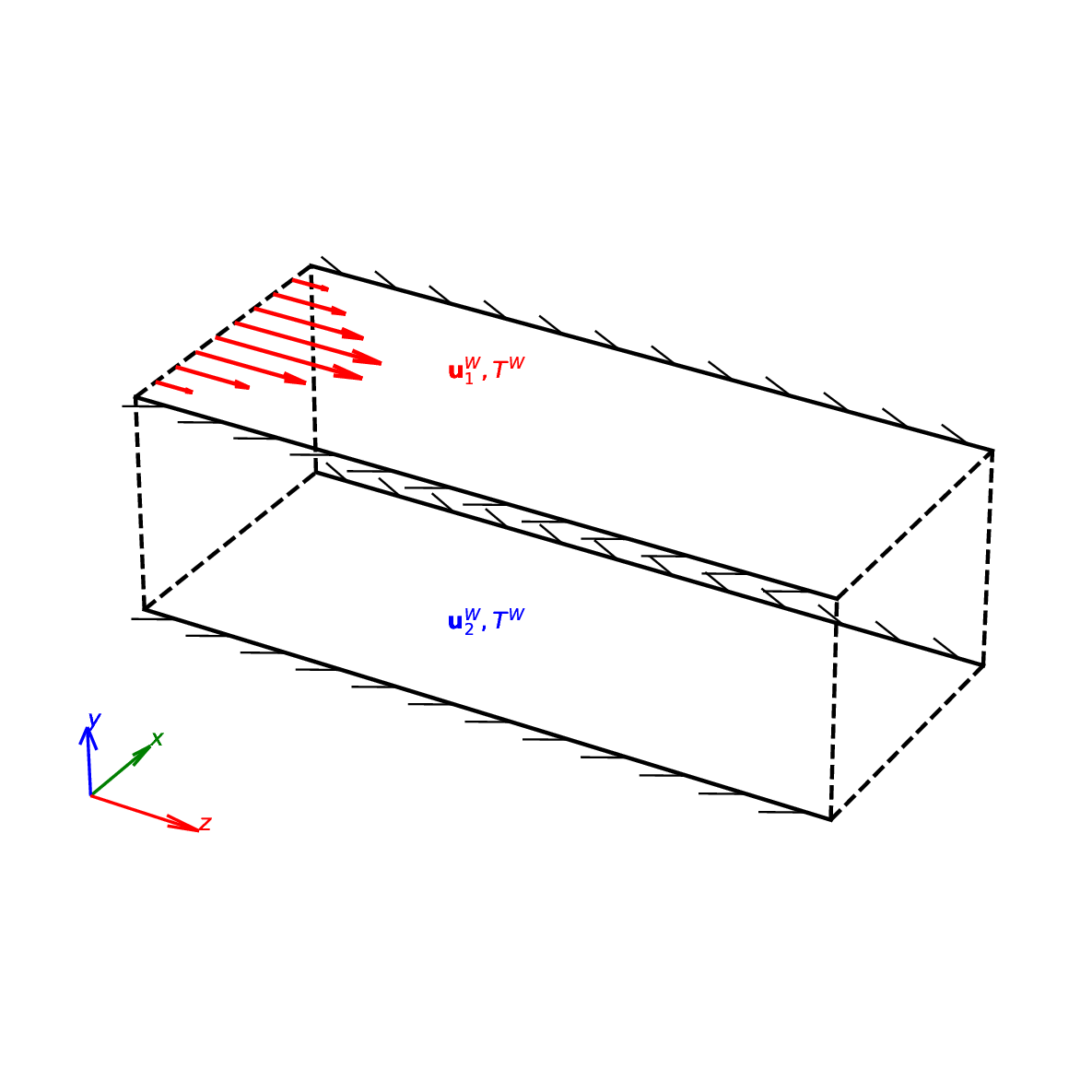}
    \caption{(2D rectangular duct flow in Sec. \ref{sec:inout}) The layout of the 2D rectangular duct flow problem.}
    \label{fig:duct_setting}
\end{figure}
In this section, the dynamics of the 2D rectangular duct flow are studied. The duct extends infinitely along the $z$-axis, and the computational domain is defined as the rectangular cross-section $\bx \in [0,1] \times [0,1]$ in the $xy$-plane. All four sides of the duct are maintained at a temperature $T^W$. Three of the duct sides are stationary with velocity $\bu^W_2 = \bz$, while the top side at $y = 1$ moves along the $z$-axis with a velocity $\bu^W_1$. In the simulation, we set 
\begin{equation}
\label{eq:duct_ini}
    \bu^W_1 = \bu_{\text{max}}\left(\sin(\pi x)\right),
\end{equation}
where $\bu_{\text{max}} = (0,0,0.5)$. Note that the velocity is continuous at the two upper corners, which can alleviate the effect of singularity on the computational accuracy. Fig. \ref{fig:duct_setting} shows the specific layout of this problem. Driven by the top plate motion, the flow eventually reaches a steady state as time approaches infinity. This setup is designed to model the dynamics of the rarefied gas flow within the duct, emphasizing the influence of a moving boundary in a controlled two-dimensional environment. The simulation requires dimension-reduced distribution functions $g$ and $h$ as well as $s_3$ since there is no symmetry in the $z$-axis, with the corresponding boundary conditions detailed as 
\begin{equation}
    \begin{aligned}
        &\displaystyle g^W(\rho^W,\boldsymbol{u}^W,T^W) = \frac{\rho^W}{{2\pi T^W}}\exp{\left(-\frac{(v_1)^2+(v_2)^2}{2T^W}\right)},\\
        &\displaystyle h^W(\rho^W,\boldsymbol{u}^W,T^W) = \frac{|\boldsymbol{u}^W|^2+T^W}{2}g^W,\\
        &\displaystyle s^W_3(\rho^W,\boldsymbol{u}^W,T^W) = (\boldsymbol{u}^W\cdot \boldsymbol{e}_3) g^W.
    \end{aligned}
\end{equation}

\begin{figure}[!hptb]
\centering
\subfigure[$\rho(\Kn=0.1)$]{
            \includegraphics[width = 0.3\textwidth]{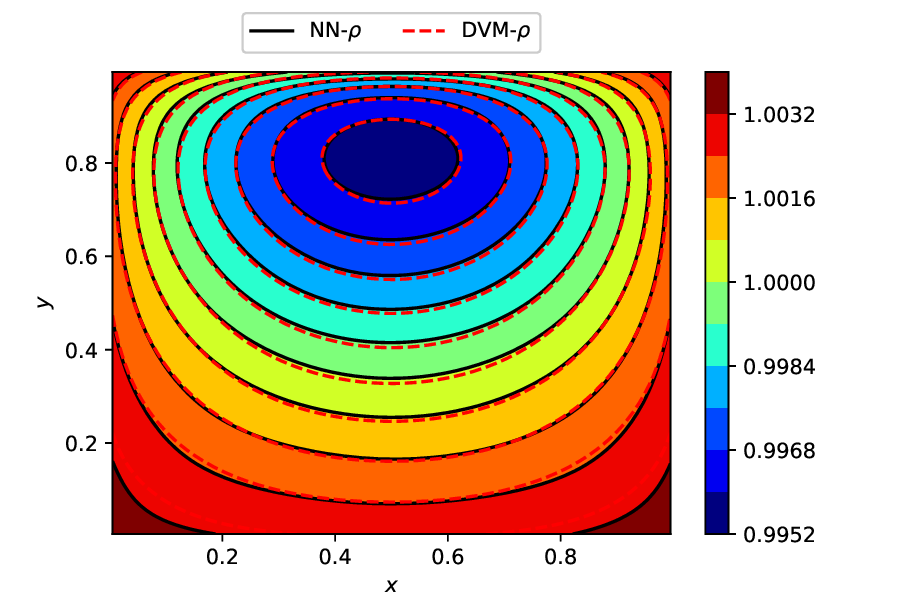}
}
\subfigure[$u_3(\Kn=0.1)$]{               
            \includegraphics[width = 0.3\textwidth]{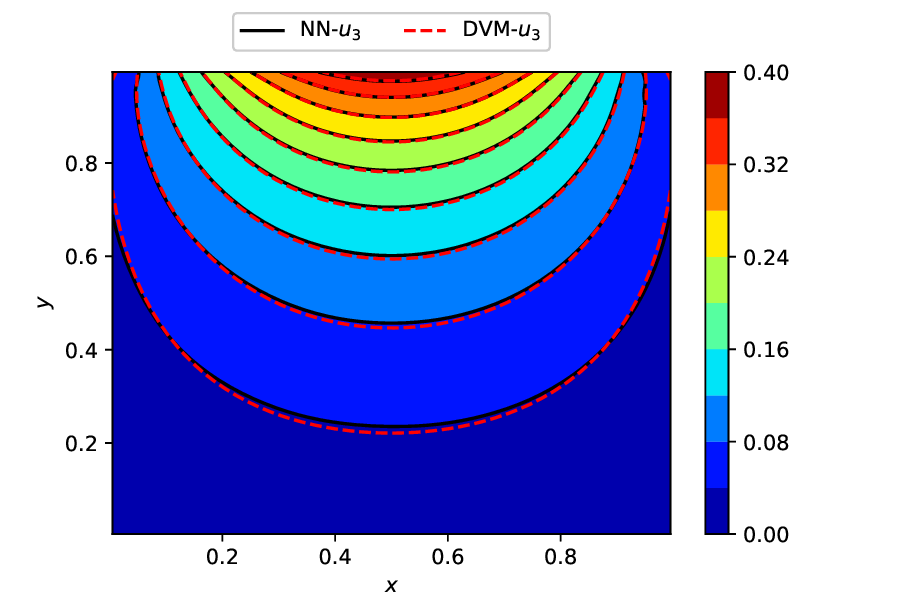}
}       
\subfigure[$T(\Kn=0.1)$]{
            \includegraphics[width = 0.3\textwidth]{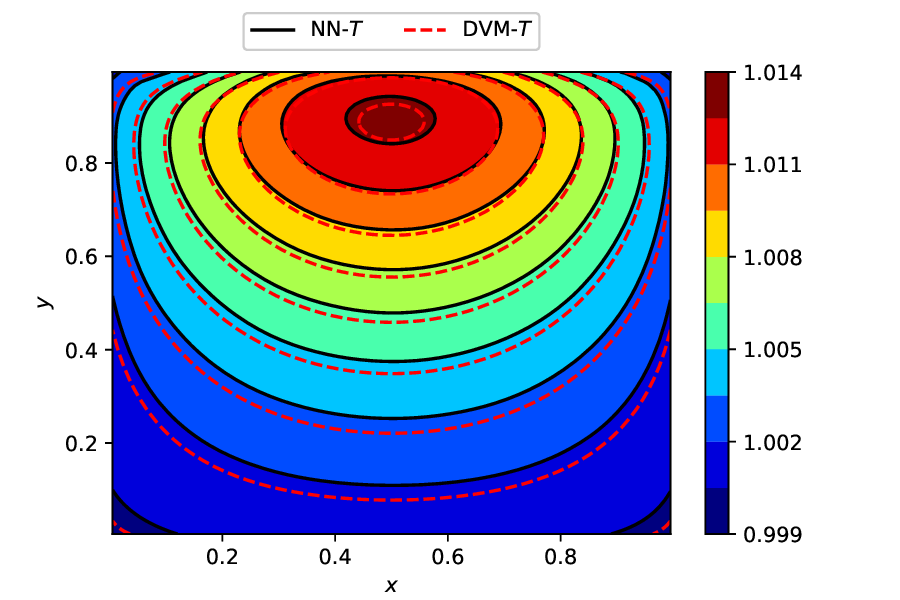}
}

\caption{(2D rectangular duct flow in Sec. \ref{sec:cavity}) Numerical solution of the density $\rho$, macroscopic velocity $u_3$, and the temperature $T$ at steady state for $\Kn = 0.1$. Here, the colored contour represents the numerical solution obtained by DRNS, while the red dashed line represents the reference solution obtained by DVM.}
    \label{fig:cavity_cdot5_kndot1}
\end{figure}

In the simulation, the spatial points are randomly selected with $N_{\text{PDE}}=2000$ in $\bx\in (0,1) \times (0,1)$ and each boundary is discretized with
$N_{\text{BC}}=300$ random points. For this 2D problem, the total training step is $20, 000$, and the activation function we used is $\sigma(x) = \sin(x)$ which demonstrated improved performance for this particular flow configuration. We first set $\Kn = 0.1$, and the numerical results of the density $\rho$, the macroscopic velocity in $z-$axis $u_3$, and the temperature $T$ at steady state are presented in Fig. \ref{fig:cavity_cdot5_kndot1}, where the reference solution obtained by DVM is also provided. It shows that for $\rho$ and $u_3$, the numerical solution by DRNS matches well with that of the reference solution, while there exists a small difference between the numerical and reference solution. But the relative error is less than $1\%$. We deduce that it is due to the non-linearity of the BGK equation when $\Kn$ is small.

\begin{figure}[!hptb]
\centering
\subfigure[$\rho(\Kn=1)$]{
            \includegraphics[width = 0.3\textwidth]{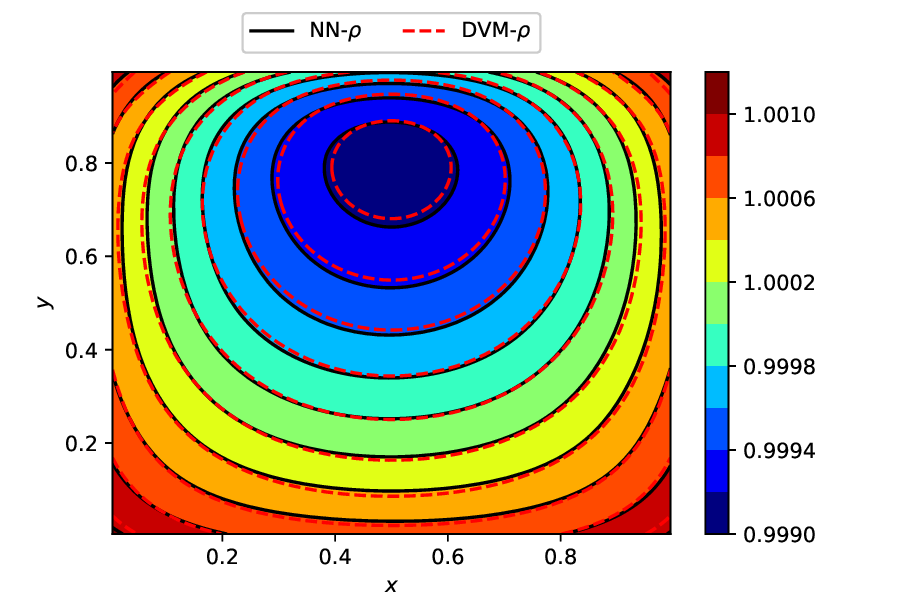}
}
\subfigure[$u_3(\Kn=1)$]{               
            \includegraphics[width = 0.3\textwidth]{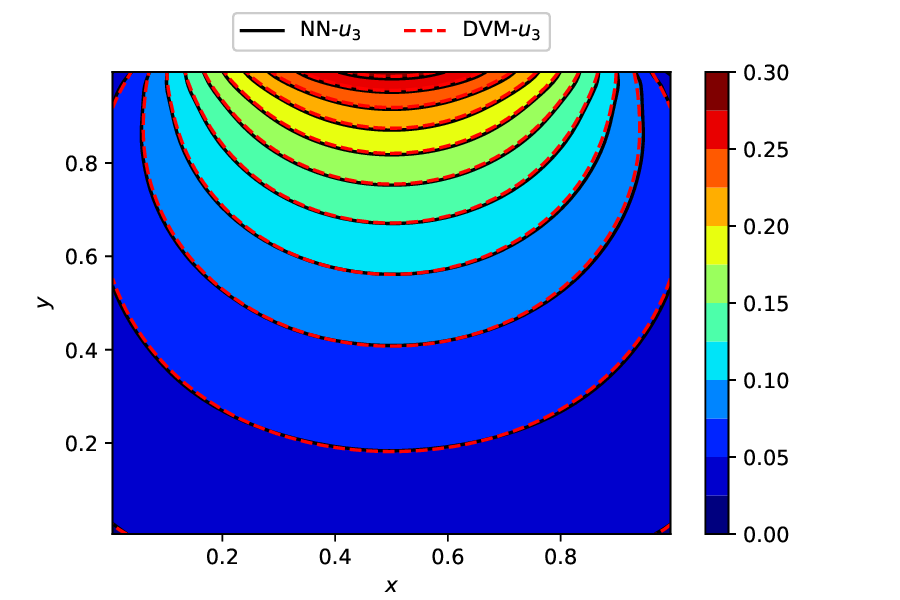}
}       
\subfigure[$T(\Kn=1)$]{
            \includegraphics[width = 0.3\textwidth]{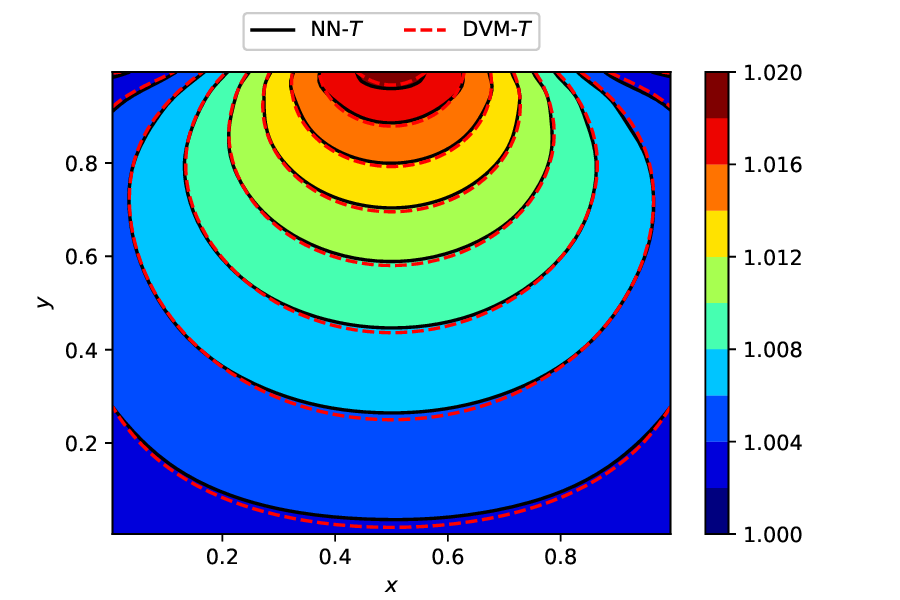}
}

\caption{(2D rectangular duct flow in Sec. \ref{sec:cavity}) Numerical solution of the density $\rho$, macroscopic velocity $u_3$, and the temperature $T$  at steady state for $\Kn = 1$. Here, the colored contour represents the numerical solution obtained by DRNS, while the red dashed line represents the reference solution obtained by DVM.}
\label{fig:cavity_cdot5_kn1}
\end{figure}

\begin{figure}[!hptb]
\centering
\subfigure[$\rho(\Kn=2.5)$]{
            \includegraphics[width = 0.3\textwidth]{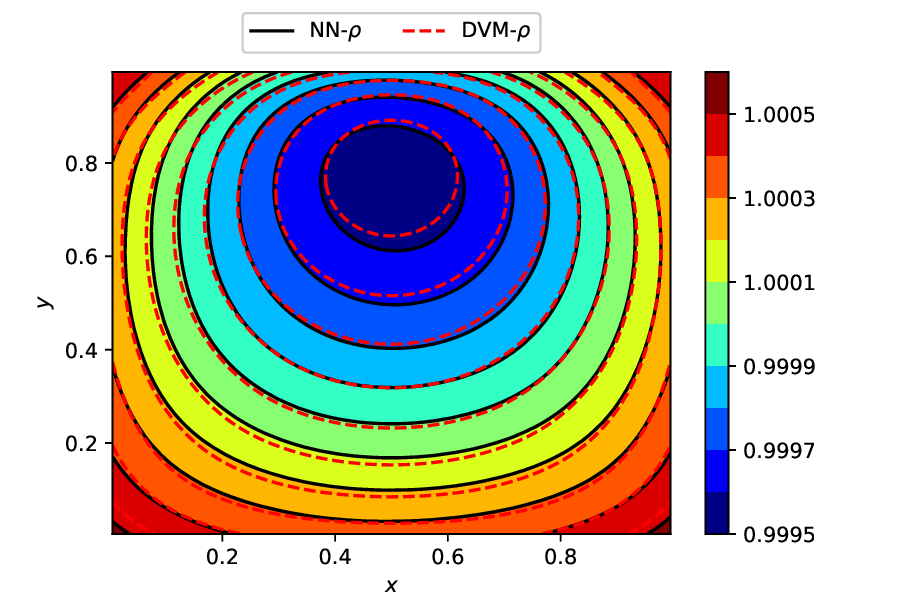}
}
\subfigure[$u_3(\Kn=2.5)$]{               
            \includegraphics[width = 0.3\textwidth]{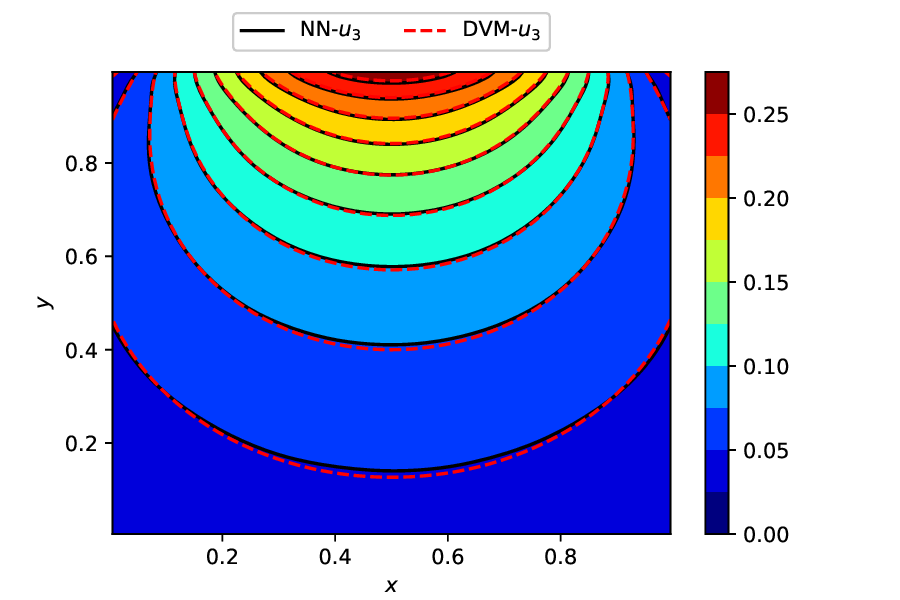}
}       
\subfigure[$T(\Kn=2.5)$]{
            \includegraphics[width = 0.3\textwidth]{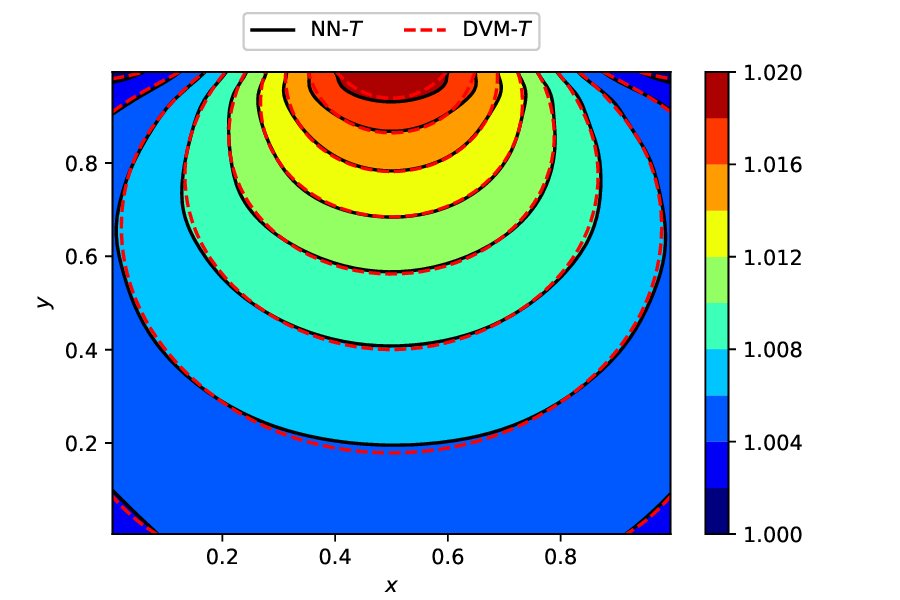}
}

\caption{(2D rectangular duct flow in Sec. \ref{sec:cavity}) Numerical solution of the density $\rho$, macroscopic velocity $u_3$, and the temperature $T$  at steady state for $\Kn = 2.5$. Here, the colored contour represents the numerical solution obtained by DRNS, while the red dashed line represents the reference solution obtained by DVM.}
\label{fig:cavity_cdot5_kn2dot5}
\end{figure}

Then, the Knudsen number is increased to $\Kn = 1$, and the numerical results of $\rho$, $u_3$ and $T$ with the reference solution obtained by DVM are provided in Fig. \ref{fig:cavity_cdot5_kn1}. With the same network and training setting, we find that the numerical solution and the reference solution are on top of each other. Finally, we increase $\Kn$ to $2.5$, which is far from the equilibrium. In this case, with the same network and training setting, the numerical solution of $u_3$ and $T$ and the reference solution by DVM are well correlated with each other, while there is a little discrepancy for the density $\rho$, as shown in Fig. \ref{fig:cavity_cdot5_kn2dot5}. This also indicates the high efficiency of this DRNS method.

\subsection{2D in-out flow}
\label{sec:inout}
\begin{figure}[!hptb]
\centering
\includegraphics[width = 0.25\textwidth]{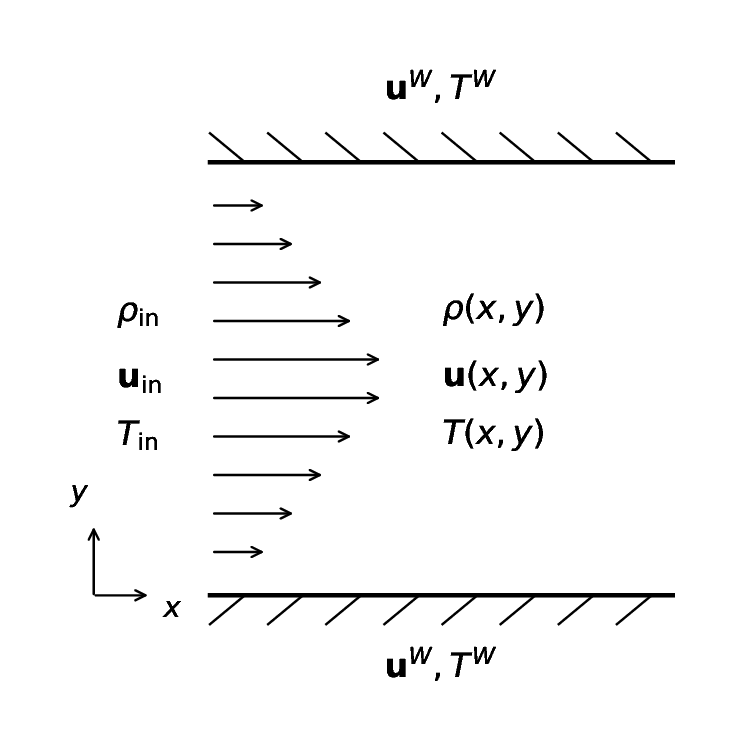}
    \caption{(2D in-out flow in Sec. \ref{sec:inout}) The layout of the 2D in-out flow problem.}
    \label{fig:inout_setting}
\end{figure}

In this section, the 2D inflow and outflow dynamics between two parallel plates is studied. The computational domain is the rectangular area $\bx \in [0,1] \times [0,1]$ between the two plates. The initial condition follows a Maxwellian distribution defined by the macroscopic variables:
\begin{equation}
    \begin{aligned}
        &\rho(x,y) = 1,\\
        &\bu(x,y) = \bu_{\max} \left(1-\tanh(10x)\right) \sin(\pi y),\\
        &T(x,y) = 1.
    \end{aligned}
\end{equation}
where $\bu_{\max} = (0.5,0,0)$. The solid boundary condition in Sec. \ref{sec:boun} is set at $y=0$ and $y=1$ with temperature $T^W = 1$ and $\bu^W = (0,0,0)$. The inflow boundary at $x=0$ follows a Maxwellian distribution with $\rho_{\rm in} = 1$, $\bu_{\rm in} = \bu_{\max} \sin(\pi y)$ and $T_{\rm in} = 1$, while the boundary at $x=1$ is designated as an outflow boundary. This setting ensures that the boundary conditions match the initial values, maintaining velocity continuity at the corners, and Fig. \ref{fig:inout_setting} shows the specific layout of this problem. The simulation requires dimension-reduced distribution functions $g$ and $h$, with the corresponding boundary conditions detailed below:
\begin{equation}
    \begin{aligned}
        &\displaystyle g^W(\rho^W,\boldsymbol{u}^W,T^W) = \frac{\rho^W}{{2\pi T^W}}\exp\left({-\frac{(\boldsymbol{u}^W\cdot \boldsymbol{e}_1-v_1)^2+(\boldsymbol{u}^W\cdot \boldsymbol{e}_2-v_2)^2}{2T^W}}\right),\\
        &\displaystyle h^W(\rho^W,\boldsymbol{u}^W,T^W) = \frac{|\boldsymbol{u}^W|^2+T^W}{2}g^W.
    \end{aligned}
\end{equation}

\begin{figure}[!hptb]
\centering
\subfigure[$\rho (\Kn = 0.1)$]{
            \includegraphics[width = 0.3\textwidth]{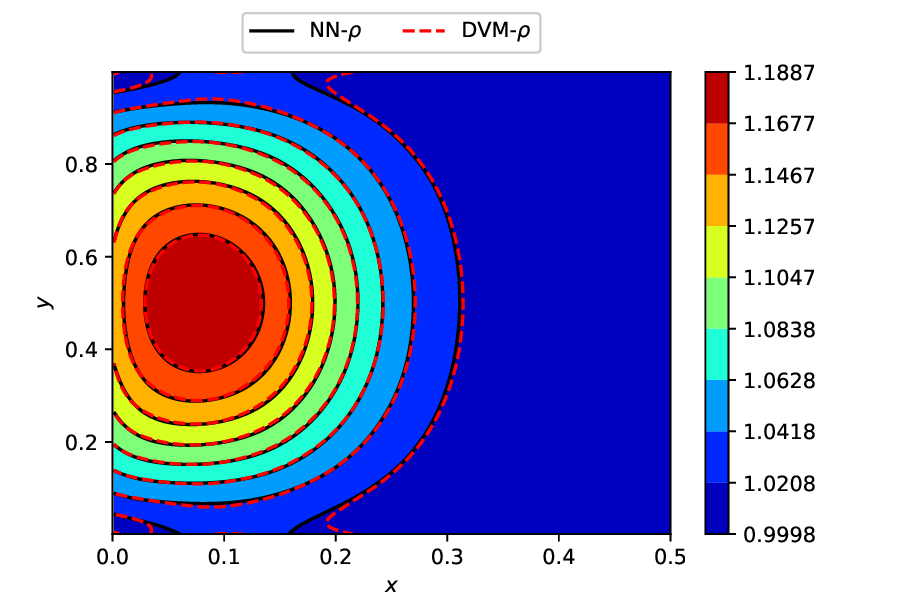}
} \hfill
\subfigure[$u_1 (\Kn = 0.1)$]{               
            \includegraphics[width = 0.3\textwidth]{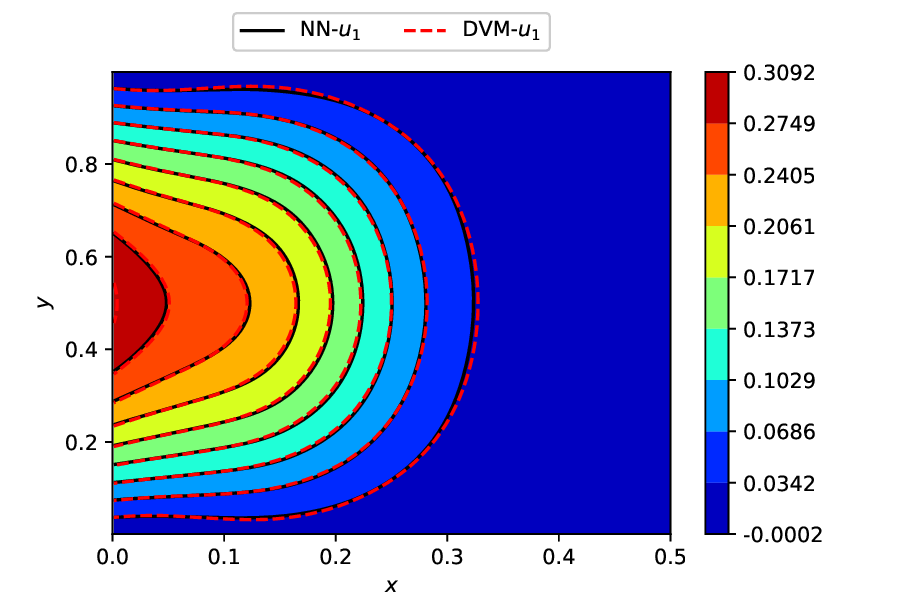}
} \hfill
\subfigure[$T (\Kn = 0.1)$]{
            \includegraphics[width = 0.3\textwidth]{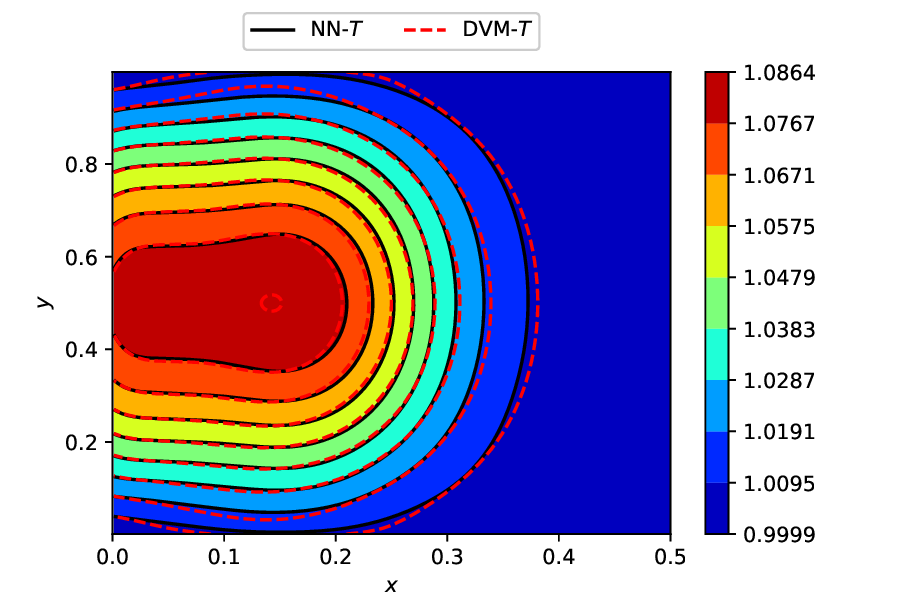}
}

\caption{(2D in-out flow in Sec. \ref{sec:inout}) Numerical solution of the density $\rho$, macroscopic velocity $u_1$, and the temperature $T$ at $t = 0.1$ for $\Kn = 0.1$. Here, the colored contour represents the numerical solution obtained by DRNS, while the red dashed line indicates the reference solution obtained by DVM. }
\label{fig:inoutflow_Cdot5_kndot1}
\end{figure}

\begin{figure}[!hptb]
\centering
\subfigure[$\rho (\Kn = 1.0)$]{
            \includegraphics[width = 0.3\textwidth]{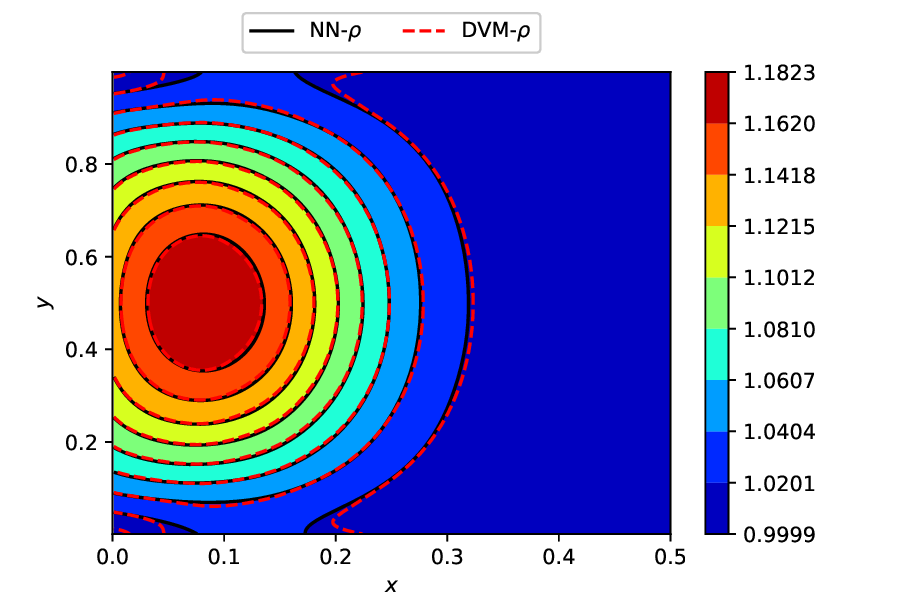}
}
\subfigure[$u_1 (\Kn = 1.0)$]{               
            \includegraphics[width = 0.3\textwidth]{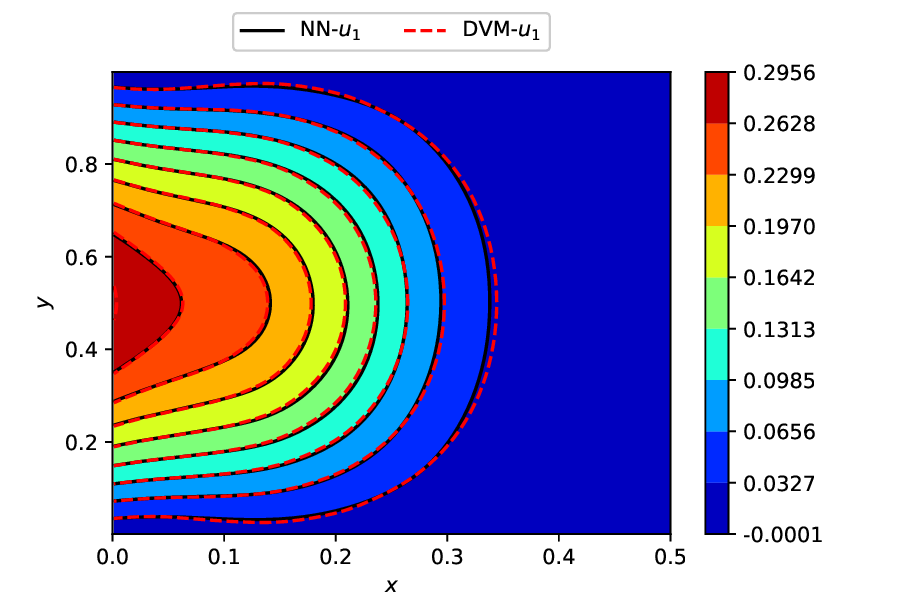}
}       
\subfigure[$T (\Kn = 1.0)$]{
            \includegraphics[width = 0.3\textwidth]{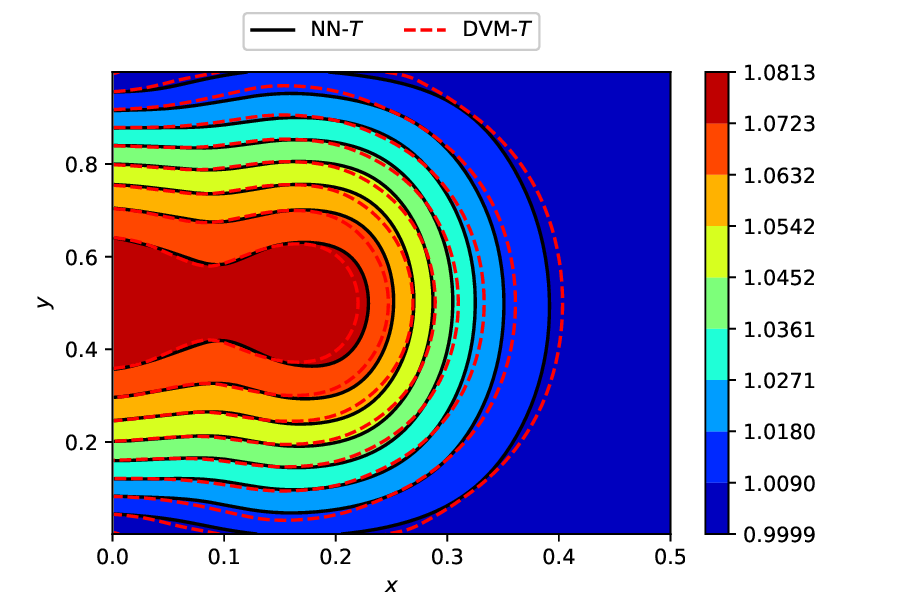}
}

	\caption{(2D in-out flow in Sec. \ref{sec:inout}) Numerical solution of the density $\rho$, macroscopic velocity $u_1$, and the temperature $T$ at $t = 0.1$ for $\Kn = 1.0$. Here, the colored contour represents the numerical solution obtained by DRNS, while the red dashed line indicates the reference solution obtained by DVM.}
    \label{fig:inoutflow_Cdot5_kn1}
\end{figure}

\begin{figure}[!hptb]
\centering
\subfigure[$\rho (\Kn = 2.5)$]{
\includegraphics[width = 0.3\textwidth]{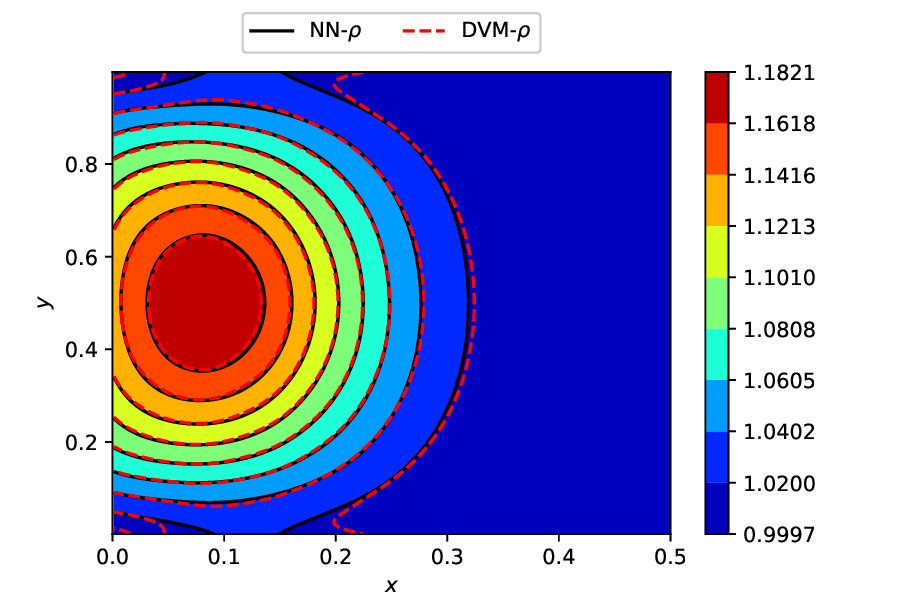}
}
\subfigure[$u_1 (\Kn = 2.5)$]{               
            \includegraphics[width = 0.3\textwidth]{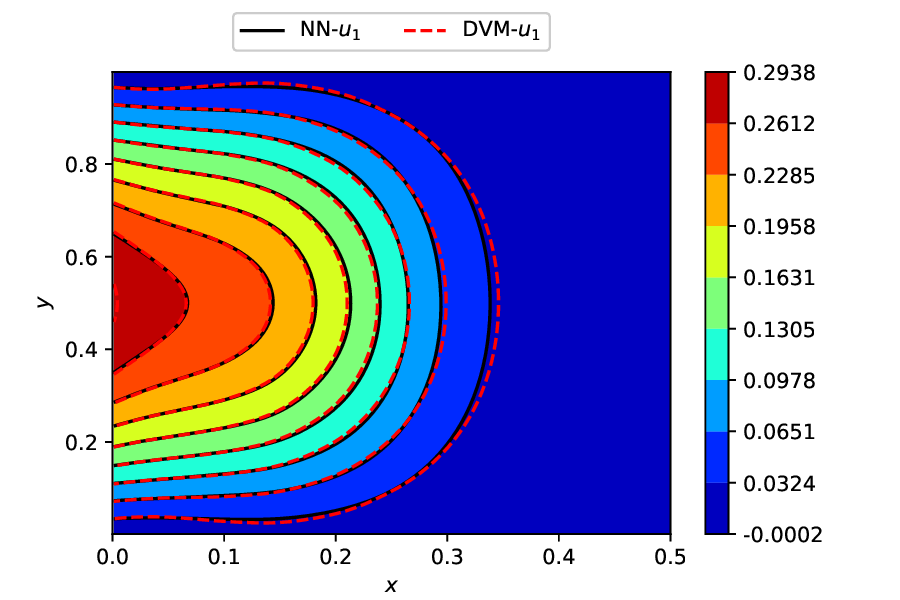}
}       
\subfigure[$T (\Kn = 2.5)$]{
            \includegraphics[width = 0.3\textwidth]{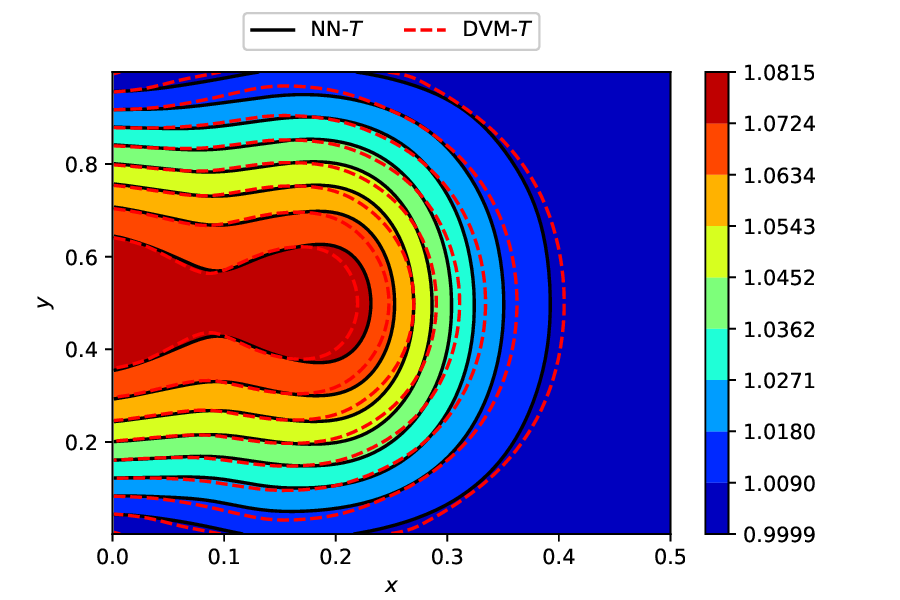}
}

\caption{(2D in-out flow in Sec. \ref{sec:inout}) Numerical solution of the density $\rho$, macroscopic velocity $u_1$, and the temperature $T$ at $t = 0.1$ for $\Kn = 2.5$. Here, the colored contour represents the numerical solution obtained by DRNS, while the red dashed line indicates the reference solution obtained by DVM.}
    \label{fig:inoutflow_Cdot5_kn2dot5}
\end{figure}

In the simulation, this time-dependent simulation spans $t\in[0,0.1]$. The number of the sampling data is set as $N_{\rm IC} = 1000$, $N_{\rm BC} = 2000$, and $N_{\rm PDE} = 3000$. This means that we randomly sample $1000$ points in $\bx \in (0, 1)^2$ at $t = 0$ for the initial condition, and $500$ points randomly sampled on each of the four boundaries during the time region $t \in (0, 0.1]$, and $3000$ randomly sampled points in the spatial space and time region $\bx \times t \in (0, 1)^2 \times (0, 0.1]$. For this 2D problem, the activation function we used is $\sigma(x) = \tanh(x)$ and the total training step is $20, 000$. We first set $\Kn = 0.1$, and the numerical results of the density $\rho$, the macroscopic velocity in $x-$axis $u_1$, and the temperature $T$ at $t = 0.1$ are shown in Fig. \ref{fig:inoutflow_Cdot5_kndot1}, where the reference solution obtained by DVM is also provided. There is an accumulation of the density $\rho$ in the middle of the region, while the macroscopic velocity $u_1$ and temperature $T$ are moving forward with the highest value along $y = 0.5$, which is also consistent with the initial condition. Fig. \ref{fig:inoutflow_Cdot5_kndot1} also indicates that the numerical solution matches well with the reference solution. Then, the Knudsen number is increased to $\Kn = 1$ and $2.5$. The numerical results at $t = 0.1$ are illustrated in Fig. \ref{fig:inoutflow_Cdot5_kn1} and \ref{fig:inoutflow_Cdot5_kn2dot5} with the reference solution obtained by DVM presented. With the same network and training setting, we find that the numerical solution and the reference solution are on top of each other, indicating the high efficiency of this DRNS method.

\section{Conclusion}
\label{sec:conclusion}
The network-based method is utilized to solve the Boltzmann-BGK equation for the microscopic flow problems. Neural representation for the distribution function is proposed, which is a high-quality ansatz to the dimensional-reduced model of the BGK equation. A specially designed network is designed to handle the Maxwell boundary in the microscopic flow problem, effectively reducing the complexity of the network parameters. The loss function including the initial and boundary condition, the residual of the PDE is proposed, resulting in significant improvement in the network's approximating efficiency. The accuracy and efficiency of this dimension-reduction neural representation method are validated through several classical benchmark problems of 1D and 2D cases.

\section*{Acknowledgments}
We thank Prof. Bin Dong from Peking University and Prof. Jun Zhang from Beihang University for their valuable suggestions. This work of Y. Wang is partially supported by the National Natural Science Foundation of China (Grant No. 12171026, U2230402 and 12031013), and President Foundation of China Academy of Engineering Physics (YZJJZQ2022017). 
%\input{article_appendix.tex}

% Bibliography
% \bibliographystyle{abbrv}  
\bibliographystyle{plain}
\bibliography{article}  

\begin{thebibliography}{10}

\bibitem{bhatnagar1954model}
P.~Bhatnagar, E.~Gross, and M.~Krook.
\newblock A model for collision processes in gases. {I}. small amplitude
  processes in charged and neutral one-component systems.
\newblock {\em Phys. Rev.}, 94(3):511--525, 1954.

\bibitem{bird1994molecular}
G.~Bird.
\newblock {\em Molecular {Gas} {Dynamics} and the {Direct} {Simulation} of
  {Gas} {Flows}}.
\newblock Oxford {Engineering} {Science} {Series}. Oxford University Press,
  1994.

\bibitem{broadwell1964study}
J.E. Broadwell.
\newblock Study of rarefied shear flow by the discrete velocity method.
\newblock {\em J. Fluid Mech.}, 19(3):401--414, 1964.

\bibitem{buet1996discrete}
C.~Buet.
\newblock {A discrete-velocity scheme for the Boltzmann operator of rarefied
  gas dynamics}.
\newblock {\em Transp. Theory Stat. Phys.}, 25(1):33--60, 1996.

\bibitem{chu1965kinetictheoretic}
C.~K. Chu.
\newblock Kinetic-{Theoretic} {Description} of the {Formation} of a {Shock}
  {Wave}.
\newblock {\em Phys. Fluids}, 8(1):12, 1965.

\bibitem{gamba2017fast}
I.~Gamba, J.~Haack, C.~Hauck, and J.~Hu.
\newblock A fast spectral method for the {Boltzmann} collision operator with
  general collision kernels.
\newblock {\em SIAM J. Sci. Comput.}, 39(4):B658--B674, 2017.

\bibitem{gamba2019micro}
I.~Gamba, S.~Jin, and L.~Liu.
\newblock Micro-macro decomposition based asymptotic-preserving numerical
  schemes and numerical moments conservation for collisional nonlinear kinetic
  equations.
\newblock {\em J. Comput. Phys.}, 382:264--290, 2019.

\bibitem{ganjaei2009new}
A.~Ganjaei and S.~Nourazar.
\newblock A new algorithm for the simulation of the {Boltzmann} equation using
  the direct simulation monte-carlo method.
\newblock {\em J. Mech. Sci. Technol.}, 23:2861--2870, 2009.

\bibitem{grad1949kinetic}
H.~Grad.
\newblock On the kinetic theory of rarefied gases.
\newblock {\em Commun. Pure Appl. Math.}, 2(4):331--407, 1949.

\bibitem{han2019uniformly}
J.~Han, C.~Ma, Z.~Ma, and W.~E.
\newblock Uniformly accurate machine learning based hydrodynamic models for
  kinetic equations.
\newblock {\em Proc. Natl. Acad. Sci. U.S.A.}, 116(44):21983--21991, 2019.

\bibitem{holloway2021acceleration}
I.~Holloway, A.~Wood, and A.~Alekseenko.
\newblock Acceleration of {B}oltzmann collision integral calculation using
  machine learning.
\newblock {\em Mathematics}, 9(12):1384, 2021.

\bibitem{hu2017asymptotic}
J.~Hu, S.~Jin, and Q.~Li.
\newblock Asymptotic-preserving schemes for multiscale hyperbolic and kinetic
  equations.
\newblock In {\em Handb. Numer. Anal.}, volume~18, pages 103--129. Elsevier,
  2017.

\bibitem{hu2020burnett}
Z.~Hu and Z.~Cai.
\newblock Burnett spectral method for high-speed rarefied gas flows.
\newblock {\em SIAM J. Sci. Comput.}, 42(5):B1193--B1226, 2020.

\bibitem{hu2020numerical}
Z.~Hu, Z.~Cai, and Y.~Wang.
\newblock Numerical simulation of microflows using {Hermite} spectral methods.
\newblock {\em SIAM J. Sci. Comput.}, 42(1):B105--B134, 2020.

\bibitem{huang2023machine}
J.~Huang, Y.~Cheng, A.J. Christlieb, L.~F. Roberts, and W.~Yong.
\newblock Machine learning moment closure models for the radiative transfer
  equation ii: {Enforcing} global hyperbolicity in gradient-based closures.
\newblock {\em Multiscale Model. Simul.}, 21(2):489--512, 2023.

\bibitem{huang2021solving}
X.~Huang, H.~Liu, B.~Shi, Z.~Wang, K.~Yang, Y.~Li, B.~Weng, M.~Wang, H.~Chu,
  J.~Zhou, F.~Yu, B.~Hua, L.~Chen, and B.~Dong.
\newblock Solving partial differential equations with point source based on
  physics-informed neural networks.
\newblock {\em arXiv:2111.01394}, 2021.

\bibitem{jin2023asymptotic}
S.~Jin, Z.~Ma, and K.~Wu.
\newblock Asymptotic-preserving neural networks for multiscale time-dependent
  linear transport equations.
\newblock {\em J. Sci. Comput.}, 94(3):57, 2023.

\bibitem{jin2024asymptotic}
S.~Jin, Z.~Ma, and K.~Wu.
\newblock Asymptotic-preserving neural networks for multiscale kinetic
  equations.
\newblock {\em Commun. Comput. Phys.}, 35(3):693--723, 2024.

\bibitem{jin2010micromacro}
S.~Jin and Y.~Shi.
\newblock A micro-macro decomposition-based asymptotic-preserving scheme for
  the multispecies {Boltzmann} equation.
\newblock {\em SIAM J. Sci. Comput.}, 31(6):4580--4606, 2010.

\bibitem{kingma2014adam}
D.~Kingma and J.~Ba.
\newblock Adam: {A} method for stochastic optimization.
\newblock {\em arXiv: 1412.6980}, 2014.

\bibitem{lecun2002efficient}
Y.~LeCun, L.~Bottou, G.~Orr, and K.~M{\"u}ller.
\newblock Efficient backprop.
\newblock In {\em Neur. Netw.: Tricks Trad.}, pages 9--50. Springer, 2002.

\bibitem{li2023hermite}
R.~Li, Y.~Lu, and Y.~Wang.
\newblock Hermite spectral method for the inelastic {Boltzmann} equation.
\newblock {\em Phys. Fluids}, 35(10), 2023.

\bibitem{li2022hermite}
R.~Li, Y.~Lu, Y.~Wang, and H.~Xu.
\newblock Hermite spectral method for multi-species {Boltzmann} equation.
\newblock {\em J. Comput. Phys.}, 471:111650, 2022.

\bibitem{li2023learning}
Z.~Li, B.~Dong, and Y.~Wang.
\newblock Learning invariance preserving moment closure model for
  {Boltzmann-BGK} equation.
\newblock {\em Commun. Math. Stat.}, 11(1):59--101, 2023.

\bibitem{li2024solving}
Z.~Li, Y.~Wang, H.~Liu, Z.~Wang, and B.~Dong.
\newblock Solving the {Boltzmann} equation with a neural sparse representation.
\newblock {\em SIAM J. Sci. Comput.}, 46(2):C186--C215, 2024.

\bibitem{lin2024monte}
Q.~Lin, C.~Zhang, X.~Meng, and Z.~Guo.
\newblock Monte {Carlo} {Physics-informed} neural networks for multiscale heat
  conduction via phonon boltzmann transport equation.
\newblock {\em arXiv preprint arXiv:2408.10965}, 2024.

\bibitem{liu2020unified}
C.~Liu and K.~Xu.
\newblock A unified gas-kinetic scheme for micro flow simulation based on
  linearized kinetic equation.
\newblock {\em Adv. Aerodyn.}, 2(1):21, 2020.

\bibitem{liu2021unified}
C.~Liu and K.~Xu.
\newblock Unified gas-kinetic wave-particle methods {IV}: multi-species gas
  mixture and plasma transport.
\newblock {\em Adv. Aerodyn.}, 3:1--31, 2021.

\bibitem{loshchilov2017sgdr}
I.~Loshchilov and F.~Hutter.
\newblock {SGDR}: {Stochastic} gradient descent with warm restarts.
\newblock {\em arXiv: 1608.03983}, 2016.

\bibitem{lou2021physicsinformed}
Q.~Lou, X.~Meng, and G.~Karniadakis.
\newblock Physics-informed neural networks for solving forward and inverse flow
  problems via the {Boltzmann}-{BGK} formulation.
\newblock {\em J. Comput. Phys.}, 447:110676, 2021.

\bibitem{maxwell1878iii}
J.C. Maxwell.
\newblock On stresses in rarefied gases arising from inequalities of
  temperature.
\newblock {\em Proc. R. Soc. Lond.}, 27(185-189):304--308, 1878.

\bibitem{miller2022neural}
S.T. Miller, N.V. Roberts, S.D. Bond, and E.C. Cyr.
\newblock Neural-network based collision operators for the {Boltzmann}
  equation.
\newblock {\em J. Comput. Phys.}, 470:111541, 2022.

\bibitem{mouhot2006fast}
C.~Mouhot and L.~Pareschi.
\newblock Fast algorithms for computing the {Boltzmann} collision operator.
\newblock {\em Math. Comput.}, 75(256):1833--1852, 2006.

\bibitem{oran1998direct}
E.S. Oran, C.K. Oh, and B.Z. Cybyk.
\newblock Direct simulation monte carlo: recent advances and applications.
\newblock {\em Annu. Rev. Fluid Mech.}, 30(1):403--441, 1998.

\bibitem{raissi2019physics}
M.~Raissi, P.~Perdikaris, and G.~Karniadakis.
\newblock Physics-informed neural networks: A deep learning framework for
  solving forward and inverse problems involving nonlinear partial differential
  equations.
\newblock {\em J. Comput. Phys.}, 378:686--707, 2019.

\bibitem{schotthofer2022neural}
S.~Schotthöfer, T.~Xiao, M.~Frank, and C.~Hauck.
\newblock Neural network-based, structure-preserving entropy closures for the
  {Boltzmann} moment system.
\newblock {\em arXiv:2201.10364}, 2022.

\bibitem{shakhov1969couette}
E.M. Shakhov.
\newblock Couette problem for the generalized {Krook} equation stress-peak
  effect.
\newblock {\em Fluid Dyn.}, 4(5):9--13, 1969.

\bibitem{shen2006rarefied}
C.~Shen.
\newblock {\em Rarefied Gas Dynamics: Fundamentals, Simulations and Micro
  Flows}.
\newblock Springer Science \& Business Media, 2006.

\bibitem{sitzmann2020implicit}
V.~Sitzmann, J.~Martel, A.~Bergman, D.~Lindell, and G.~Wetzstein.
\newblock Implicit neural representations with periodic activation functions.
\newblock In {\em Adv. Neural Inf. Process Syst.}, volume~33, pages 7462--7473,
  2020.

\bibitem{struchtrup2005macroscopic}
H.~Struchtrup.
\newblock {\em Macroscopic transport equations for rarefied gas flows:
  approximation methods in kinetic theory}.
\newblock Interaction of mechanics and mathematics series. Springer, 2005.

\bibitem{wang2019approximation}
Y.~Wang and Z.~Cai.
\newblock Approximation of the {Boltzmann} collision operator based on hermite
  spectral method.
\newblock {\em J. Comput. Phys.}, 397:108815, 2019.

\bibitem{wu2013deterministic}
L.~Wu, C.~White, T.~Scanlon, J.~Reese, and Y.~Zhang.
\newblock Deterministic numerical solutions of the {Boltzmann} equation using
  the fast spectral method.
\newblock {\em J. Comput. Phys.}, 250:27--52, 2013.

\bibitem{xiao2021using}
T.~Xiao and M.~Frank.
\newblock Using neural networks to accelerate the solution of the {Boltzmann}
  equation.
\newblock {\em J. Comput. Phys.}, 443:110521, 2021.

\bibitem{xu2019frequency}
Z.~Xu, Y.~Zhang, T.~Luo, Y.~Xiao, and Z.~Ma.
\newblock Frequency principle: {Fourier} analysis sheds light on deep neural
  networks.
\newblock {\em arXiv preprint arXiv:1901.06523}, 2019.

\bibitem{yang1995rarefied}
J.~Yang and J.~Huang.
\newblock Rarefied flow computations using nonlinear model {Boltzmann}
  equations.
\newblock {\em J. Comput. Phys.}, 120(2):323--339, 1995.

\bibitem{zhang2023simulation}
L.~Zhang, W.~Ma, Q.~Lou, and J.~Zhang.
\newblock Simulation of rarefied gas flows using physics-informed neural
  network combined with discrete velocity method.
\newblock {\em Phys. Fluids}, 35(7), 2023.

\end{thebibliography}

\end{document}